\documentclass[12pt]{iopart}

\usepackage[T1]{fontenc}
\usepackage{tgadventor}
\usepackage{times}
\usepackage{url}
\usepackage{amsmath}
\usepackage{amssymb}
\usepackage{amsthm}
\usepackage{float}
\usepackage{tikz}
\usetikzlibrary{arrows, snakes, backgrounds}
\usepackage{graphicx}
\usepackage{caption}
\usepackage{subcaption}
\usepackage{pgfmath}
\usepackage[boxed,resetcount]{algorithm2e}
\newcommand{\R}[1]{\mathbb{R}^{#1}}
\usetikzlibrary{positioning,arrows.meta,quotes}
\newcommand{\HT}{$\mathcal{H}^{2}$~}
\newcommand*{\hermconj}{^{\mathsf{H}}}
\newcommand*{\conj}[1]{\overline{#1}}
\newtheorem{remark}{Remark}
\definecolor{darkred}{RGB}{173,0,43}
\definecolor{dgreen}{RGB}{0, 115, 85}
\definecolor{lgreen}{RGB}{217, 231, 225}
\definecolor{lblue}{RGB}{0, 115, 85} 

\begin{document}

\title[Multidimensional deconvolution with shared bases]{Multidimensional deconvolution with shared bases}

\author{Daria Sushnikova, Matteo Ravasi and David Keyes\\
Extreme Computing Research Center \\
King Abdullah University of Science and Technology}

\address{0138-B1, KAUST,
Thuwal, Saudi Arabia, zip 23955}
\ead{daria.sushnikova@gmail.com}
\vspace{10pt}
\begin{indented}
\item[] April 1, 2024 
\end{indented}

\begin{abstract}
We address the estimation of seismic wavefields by means of Multidimensional Deconvolution (MDD) for various redatuming applications. While offering more accuracy than conventional correlation-based redatuming methods, MDD faces challenges due to the ill-posed nature of the underlying inverse problem and the requirement to handle large, dense, complex-valued matrices. These obstacles have long limited the adoption of MDD in the geophysical community. Recent interest in this technology has spurred the development of new strategies to enhance the robustness of the inversion process and reduce its computational overhead. We present a novel approach that extends the concept of block low-rank approximations, usually applied to linear operators, to simultaneously compress the operator, right-hand side, and unknowns. This technique greatly alleviates the data-heavy nature of MDD. Moreover, since in 3d applications the matrices do not lend themselves to global low rank approximations, we introduce a novel \HT-like approximation. We aim to streamline MDD implementations, fostering efficiency and controlling accuracy in wavefield reconstruction. This innovation holds potential for broader applications in the geophysical domain, possibly revolutionizing the analysis of multi-dimensional seismic datasets.

\end{abstract}

%
%
%
%
%

\section{Introduction}
Seismic recordings provide crucial insights into the interaction of waves with an otherwise unknown physical medium. A so-called wave state (defined by the wave nature, medium parameters, source type, and boundary conditions) can be associated with any such recording~\cite{fokk-intr-2013}. While, in theory, one has the flexibility to choose any combination of such conditions, real-world seismic exploration faces constraints stemming from physical, economical, and environmental factors.

Several data-driven approaches have been developed under the framework of seismic interferometry (e.g.,~\cite{wape-gf-2006}) to create ideal settings, also known as virtual surveys, from physically acquired data. By defining interaction quantities and employing reciprocity theorems~\cite{fokk-intr-2013}, these approaches can leverage wavefield measurements from distinct receivers to emulate responses as if a source was placed in one of the receivers' locations. Ideally, cross-correlating the wavefield recorded by two receivers surrounded by a closed boundary of sources in a lossless medium yields the Green’s function between such receivers. This principle has been extensively used in seismic applications for a variety of tasks~\cite{li-intr2-2021,ruig-intr2-2010,drag-intr2-2009,schu-intr2-2006, baku-intr-2006,schu-intr-2004}. 

However, closed boundaries are rarely available in real-life applications. This, combined with the inevitable presence of losses in the propagating medium, can hinder accurate reconstruction, yielding a blurred Green's function from cross-correlation interferometry~\cite{wape-intr-2011}. An alternative approach that can accommodate both dissipative media and open boundaries relies on a convolution-type reciprocity theorem. ~\cite{wape-intr-2011} showed that this leads to an implicit Green's function representation; deconvolution can then be employed to retrieve such Green's function mitigating the aforementioned correlation-induced, spatio-temporal blurring effects. This method, known as seismic interferometry by a Multidimensional Deconvolution (MDD) has been initially proposed in the acoustic setting and later extended to more general vector fields, including elastic~\cite{van-elast-2011} and electromagnetic waves~\cite{wape-em-2008}.

While the adoption of MDD constitutes a significant leap towards achieving more accurate redatuming compared to traditional correlation methods, it comes with numerical and computational challenges. The former stems from the ill-posed nature of the inverse problem that one has to solve to estimate the sought after Green's function, whilst the latter is related to the arithmetic and storage complexity of the problem. MDD algorithms must in fact manipulate dense matrices for the operator, unknowns, and right-hand side, which significantly escalates the computational demands as the problem size grows, especially in 3d applications. This has resulted in a somewhat constrained uptake of MDD within the geophysical community.

However, the superior capabilities of MDD have led the geophysical community to identify strategies aimed at bolstering inversion stability and alleviating computational overhead. For instance, to enhance stability, researchers in~\cite{mina-mdd_svd-2011} apply singular-value decomposition (SVD) to compute the pseudoinverse matrix and truncate small eigenvalues as a way to naturally regularize the inverse process. More recently, ~\cite{boiero-mdd_stab-2021} proposed to introduce a regularization term to enforce solution similarity among neighboring sources and/or receivers. 

Alternatively, a time-domain formulation of MDD was proposed in~\cite{neut-imtf-2013}: by operating in the time domain, one can more naturally
regularize the inverse problem, alongside parametrizing the solution in a suitable transformed domain (e.g., curvelets) where sparsity-promoting inversion can be employed. Along similar lines, \cite{neut-targ-2017} and~\cite{luik-ssmc-2020} proposed to introduce a preconditioner to enforce reciprocity
in the sought after Green's function. Continuing,~\cite{varg-mdd-2021} suggested to combine three physics-based preconditioners, namely causality, reciprocity, and frequency–wavenumber locality within the same time-domain formulation. This resulted in improved stability also for scenarios with highly complex propagating overburdens (e.g., salt bodies). Finally, low-rank regularization has been recently proposed by~\cite{kuma-rtrans-2022,chen-qq-2023,chen-nn-2023} as an alternative way to parametrize the sought after solution. Whilst primarily focused on tackling the first challenge (i.e., stability), this approach can also reduce the size of the solution when a factorization-based form of low-rankness is considered.

From a computational standpoint, two main approaches have been investigated to date to reduce the cost of applying MDD. The former leverages the fact that the MDD equations can be formulated as a finite-sum functional~\cite{rava-stmdd-2022}; as such, the associated inverse problem can be solved by means of stochastic gradient descent algorithms, where
the gradient at each iteration is computed using a small subset of randomly selected sources. This can reduce the overhead of having to access the entire dataset at each iteration and has been shown to lead to faster and more robust convergence. The latter involves the application of block low-rank approximation techniques, such as the Tile-Low Rank (TLR) algorithm~\cite{hong-tlr-2021,ltai-smw-2023,hong-hpc-2022,rava-large-2022,rava-tlr-2022}, to the operator matrix. By doing so, one can first reduce the size of the operator, and later accelerate matrix-vector multiplications to performing this operation in two steps where the factors, instead of the full matrix, are utilized.

Our paper extends the concept of block low-rank approximation to encompass compression of both the operator, itself, and the right-hand side and unknown matrix. This is particularly advantageous given the data-heavy nature of the right-hand side and the unknowns in MDD. Furthermore, by embracing a more sophisticated \HT~\cite{hack-h2-2000,grro-fmm-1987} approximation of the dense matrices, we show that further compressibility of the wavefield can be achieved compared to the TLR approach. Through this multifaceted approach, we endeavor to pave the way for more efficient MDD implementations, ultimately expanding its potential applications within the geophysical domain.

\section{Multi-dimensional deconvolution}

\begin{figure}[H]
\begin{center}
\begin{tikzpicture}[scale=1.]
    \begin{scope}[shift={(0,0)}]
        \draw [draw=dgreen,fill=dgreen] (0, 0) rectangle (10, 8);
        \draw [draw=dgreen,fill=lgreen] (0, 4) rectangle (10, 8);
        \draw[step=2.0, dgreen,thin] (0,4) -- (10,4);
        \foreach \x in {1,...,9}{
            \pgfmathrandominteger{\a}{-240}{0}
            \ifthenelse{\x=3}
            {\filldraw[red] (\x, 8) circle (2pt)}
            {\filldraw[red] (\x, \a/240+8) circle (2pt)};
        }
        \foreach \x in {0,...,9}{
            \filldraw[lblue] (\x+0.5, 4) circle (2pt);
        }
        \foreach \x in {1,...,9}{
            \filldraw[white] (\x, 4) circle (2pt);
        }
        \node at (5,9.5) {\color{dgreen}\Large ground or water surface 1d or 2d};
        \draw[->, dgreen] (5,9.2) -- (5, 8.3);
        \node at (2.5,8.75) {\color{red}\Large sources};
        \node at (7.2,2.8) {\color{lgreen}\Large receivers};
        \node at (3.1,7.50) {\color{red}\Large $x_s$};
        \draw[->, red] (2.5,8.5) -- (2.8, 8.2);
        \draw[->, lgreen] (7,3.2) -- (7, 3.8);
        \node at (5,6) {\color{dgreen}\Large overburden};
        \node at (1,4.5) {\color{lblue}\Large $x'_0$};
        \node at (3.5,4.5) {\color{lblue}\Large $x_0$};
        \node at (5,2) {\color{lgreen}\Large target of interest};
        \draw[->, dgreen] (7.5,4.7) -- (7.5, 4.2);
        \node at (7.5,5) {\color{lblue}\Large virtual receivers};
        \draw[lgreen] (1,4) -- (2.25, 3);
        \draw[->, lgreen] (2.25,3) -- (3.4, 3.9);
        \node at (3.2,3.2) {\color{lgreen}\Large $r$?};
    \end{scope}
\end{tikzpicture}
\end{center}
\caption{Visual representation of the MDD algorithm.}
\label{fig:mdc}
\end{figure}
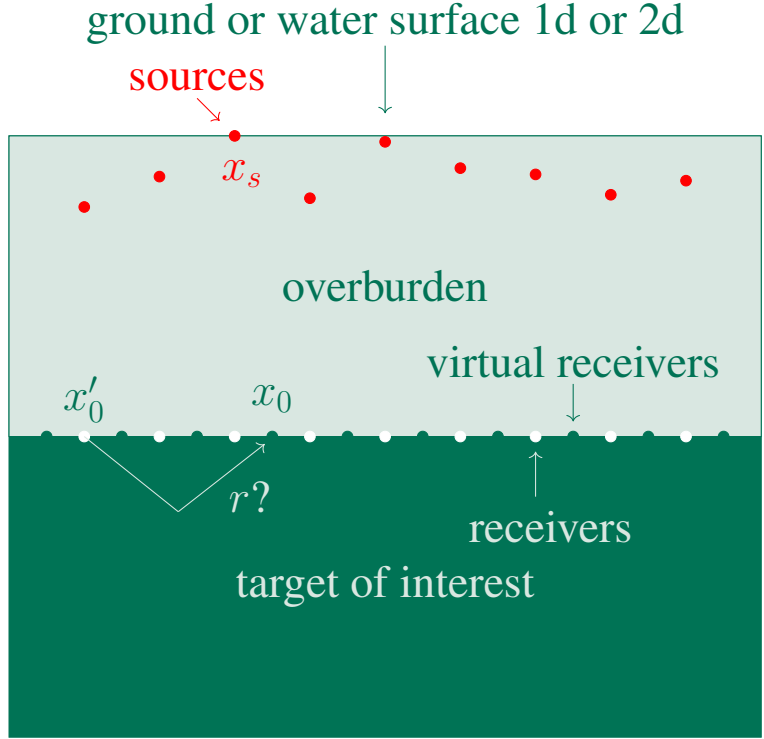

Let us consider a medium composed of two parts: the first one, usually referred to as overburden, lies in between the surface of the Earth and the receiver datum; the second one, below the receiver datum, represents the target area we are ultimately interested to image. One or more sources, denoted here as $x_s$, are placed anywhere in the overburden, whilst receivers $x_0'$ and virtual sources $x_0$ are located at the datum that separates the overburden from the target area. We assume that the wavefields $w^{-}$ from $x_s$ to $x_0$ and $w^{+}$ from $x_s$ to $x_0'$ are known (i.e., physically recorded or obtained by a pre-processing step). In this problem, we are interested to recover the wavefield $r$ that travels only below the receiver datum. Figure~\ref{fig:mdc} shows a schematic representation of the acquisition geometry associated with the MDD problem.


Following~\cite{wape-ssmc-2011}, the aforementioned wavefields are related to each other via a convolution-type representation theorem,

\begin{equation} \label{eq:mdccont}
w^{-} (x_{s}, x_{0}, t)= \mathcal{F}_t^{-1} \int w^{+}(x_s, x'_0,  f)\mathcal{F}_t \left(r(x'_0, x_0, t)\right) dx'_0, 
\end{equation}
where $\mathcal{F}_t$ and $\mathcal{F}_t^{-1}$ represent the forward and inversion Fourier transforms, whilst $t$ and $f$ refer to the time and frequency axes. In order to find a numerical solution of equation~(\ref{eq:mdccont}) for the wavefield of interest $r$, both space and time are usually discretized, leading to the following expression

\begin{equation} 
\label{eq:mdc}
w^{-} (x_{s,i}, x_{0,j}, t_k)= F_t^{-1} \left( \sum_{\alpha=1}^{N_0} w^{+}(x_{s,i}, x'_{0,\alpha}, f_p) F_t (r(x'_{0,\alpha}, x_{0,j}, t_k))\right), 
\end{equation}
where $F_t$ and $F_t^{-1}$ are the discrete forward and inversion Fourier transforms, 
$i=1\dots N_s$, $j=1,\dots, N_0$, $\alpha=1,\dots, N_0$, $k=1,\dots, N_t$, and $p=1,\dots, N_f$,
with $N_s$, $N_0$, $N_t$, $N_f$ representing the number of sources, receivers, time samples, and frequency samples, respectively.

The linear system in equation~(\ref{eq:mdc}) can be reformulated in a tensor form if we assume that
$$B_t[i, j, k] = w^{-}(x_{s, i}, x_{0, j}, t_k)$$
$$W[i,\alpha,p] = W^{+}(x_{s, i}, x'_{0,\alpha}, f_p)$$
$$X_t[\alpha,j,k] = r(x'_{0,\alpha}, x_{0, j}, t_k)$$
yielding
$$F^H W F X_t = B_t,$$
where $F$ is a Discrete Fourier transform operator applied over the time index, $W\in\R{N_s\times N_0\times N_f}, X_t\in\R{N_0\times N_0\times N_t}, B_t\in\R{N_s\times N_0\times N_t}$.

By applying the Fourier transform to the system from right and left over the time index, we obtain an equivalent system of equations in the frequency domain:

$$W X_f = B_f,$$
where $X_f\in\R{N_0\times N_0\times N_f}$ and $B_f\in\R{N_s\times N_0\times N_f}$. This tensor system is equivalent to a series of linear systems for different frequencies:
\begin{equation}
W_p X_p = B_p, \quad p=1\dots N_f,
\label{eq:fr_sys}
\end{equation}
where $W_p = W[:,:, p]$, $X_p = X_f[:,:, p]$, and $B_p = B_f[:,:, p]$.

\begin{center}
\begin{tikzpicture}[scale=0.5]
    \begin{scope}[shift={(0,0)}]
    \draw [draw=dgreen, fill=lgreen] (0, 0) rectangle (4, 4);
    \node at (2,4.5) {$N_s\times N_0$};
        \node at (2,-1) {$W_p$};
    \end{scope}
    \begin{scope}[shift={(5,0)}]
        \draw [draw=dgreen, fill=lgreen] (0, 0) rectangle (4, 4);
        \node at (2,4.5) {$N_0\times N_0$};
        \node at (5.,2) {\huge $=$};
        \node at (2,-1) {$X_p$};
    \end{scope}
    \begin{scope}[shift={(11,0)}]
        \draw [draw=dgreen, fill=lgreen] (0, 0) rectangle (4, 4);
        \node at (2,4.5) {$N_s\times N_0$};
        \node at (2,-1) {$B_p$};
    \end{scope}
\end{tikzpicture}
\end{center}

The three matrices involved in equation~\ref{eq:fr_sys} are usually dense and have comparable sizes. In order to reduce the cost of solving these $N_f$ systems of equations, one must find a suitable compressed representation of such matrices.

\section{Compression of the operator, right-hand side, and unknowns}
In this section, we introduce several approaches that benefit from the global or block-wise low-rank nature of the matrices in the system of equations~(\ref{eq:fr_sys}). These approaches span from the simplest low-rank approximation to more advanced techniques based on the \HT matrix approximation. 

\subsection{Global low-rank parametrization}
To begin, let us examine the scenario where all the matrices in equation~(\ref{eq:fr_sys}) possess the property of being globally low rank.  

\subsubsection{Approximation based on the non-symmetric low-rank parametrization of~$X_p$}
Since we expect the solution $X_p$ to be globally low rank, let us write $X_p$ in this factorized form:
$$X_p = L_pR_p,$$
where $X_p\in\R{N_0 \times N_0}$, $L_p\in\R{N_0\times r}$, $R_p\in\R{r\times N_0}$, $r\ll N_0$.
We obtain the system:
\begin{equation}
    W_pL_pR_p = B_p.
    \label{eq:xlr1}
\end{equation}

One obvious option for solving this system of equations is represented by the alternating least squares method (ALS~\cite{taul-als-1993,taul-als-1995}), where the factors $L_p$ and $R_p$ are updated in an alternating fashion. This method has been shown to have poor convergence unless additional regularization terms are imposed on the factors (see, for example,~\cite{chen-qq-2023}). Instead, we decide here to compute also a low-rank approximation for the matrix $B_p$:

$$B_p = \widehat{L}_p\widehat{R}_p,$$
where $B_p\in\R{N_s \times N_0}$, $\widehat{L}_p\in\R{N_s\times r}$, $\widehat{R}_p\in\R{r\times N_0}$, $r\ll N_0, N_s$. 
We obtain the system:
\begin{equation}
    W_pL_pR_p = \widehat{L}_p\widehat{R}_p.
    \label{eq:xlr1}
\end{equation}

\begin{center}
\begin{tikzpicture}[scale=0.5]
    \begin{scope}[shift={(0,0)}]
    \draw [draw=dgreen, fill=lgreen] (0, 0) rectangle (4, 4);
    \node at (2,4.5) {$N_s\times N_0$};
        \node at (2,-1) {$W_p$};
    \end{scope}
    \begin{scope}[shift={(5,0)}]
        \draw [draw=dgreen, fill=lgreen] (0, 0) rectangle (1, 4);
        \node at (0.5,4.5) {$N_0\times r$};
        \node at (0.5,-1) {$L_p$};
    \end{scope}
    \begin{scope}[shift={(7,0)}]
        \draw [draw=dgreen, fill=lgreen] (0, 3) rectangle (4, 4);
        \node at (2,4.5) {$r\times N_0$};
        \node at (5.,2) {\huge $=$};
        \node at (2,-1) {$R_p$};
    \end{scope}
    \begin{scope}[shift={(13,0)}]
        \draw [draw=dgreen, fill=lgreen] (0, 0) rectangle (1, 4);
        \node at (0.5,4.5) {$N_s\times r$};
        \node at (0.5,-1) {$\widehat{L}_p$};
    \end{scope}
    \begin{scope}[shift={(15,0)}]
        \draw [draw=dgreen, fill=lgreen] (0, 3) rectangle (4, 4);
        \node at (2,4.5) {$r\times N_0$};
        \node at (2,-1) {$\widehat{R}_p$};
    \end{scope}
\end{tikzpicture}
\end{center}


We then assume that
\begin{equation}
    R_p = \widehat{R}_p.
    \label{eq:xlr2}
\end{equation}

Such an assumption is based on the intuition that the solution and right-hand side are wavefields with similar features and therefore can be represented by a shared basis. The experiments in sections~\ref{subsec:ds1} and~\ref{subsec:ds2} confirm the validity of this assumption as satisfactory solutions are produced in both scenarios. Limitations and benefits of such an assumption are a subject of future research.

Using the assumption in equation~(\ref{eq:xlr2}), equation~(\ref{eq:xlr1}) becomes:

\begin{equation}
    W_pL_p\widehat{R}_p = \widehat{L}_p\widehat{R}_p.
    \label{eq:xlr3}
\end{equation}
Multiplying both sides of equation~(\ref{eq:xlr3})\ by the pseudoinverse of $\widehat{R}_p$ (i.e., $\widehat{R}_p^{\dagger}$) we obtain:

\begin{equation}
    W_pL_p = \widehat{L}_p
    \label{eq:xlr4}
\end{equation}

\begin{center}
\begin{tikzpicture}[scale=0.5]
    \begin{scope}[shift={(0,0)}]
        \draw [draw=dgreen, fill=lgreen] (0, 0) rectangle (4, 4);
        \node at (2,4.5) {$N_s\times N_0$};
        \node at (2,-1) {$W_p$};
    \end{scope}
    \begin{scope}[shift={(5,0)}]
        \draw [draw=dgreen, fill=lgreen] (0, 0) rectangle (1, 4);
        \node at (0.5,4.5) {$N_0\times r$};
        \node at (0.5,-1) {$L_p$};
        \node at (2.,2) {\huge $=$};
    \end{scope}
    \begin{scope}[shift={(8,0)}]
        \draw [draw=dgreen, fill=lgreen] (0, 0) rectangle (1, 4);
        \node at (0.5,4.5) {$N_s\times r$};
        \node at (0.5,-1) {$\widehat{L}_p$};
    \end{scope}

\end{tikzpicture}
\end{center}

We can now solve this system for $L_p$ and ultimately reconstruct our solution of interest as $X_p = L_p\widehat{R}_p$. By choosing this parametrization, we have reduced the size of the problem (and therefore the number of computations required to obtain the solution) by a factor of $\frac{N_0}{r}$. We refer to this technique as the USV method and as USV-rsp when a reciprocity preconditioner is added in the solution of equation \ref{eq:xlr4}.

However, the proposed approach presents several downsides. Notably, it does not exploit the fact that also $W_p$ may be a low-rank matrix. Moreover, the chosen parametrization for the solution (i.e., $X_p = L_p \widehat{R}_p$) ignores the fact that this matrix should be symmetric, owing to the reciprocity of acoustic wavefields. 


\subsubsection{Approximation based on the symmetric low-rank parametrization of $X_p$.}
\label{sec:xvsv}
Let us now define the truncated SVD of the matrix $B_p$ as:

$$B_p = U_p\Sigma_pV_p,$$
where $B_p\in\R{N_s \times N_0}$, $U_p\in\R{N_s\times r}$, $\Sigma_p\in\R{r\times r}$, $V_p\in\R{r\times N_0}$,  $r\ll N_0, N_s.$ 

Given that we want i) $X_p$ to be symmetric (but not Hermitian), and ii) the matrices $B_p$ and $X_p$ to share their right factor, we write $X_p$ in the following format: 
$$X_p = V_p^{\top}S_p V_p,$$
where $S_p\in\R{r\times r}$ is an unknown matrix. We obtain the system:

\begin{equation}
    W_p V^{\top}_p S_p V_p =  U_p\Sigma_pV_p.
    \label{eq:vsv1}
\end{equation}

\begin{center}
\begin{tikzpicture}[scale=0.5]
    \begin{scope}[shift={(0,0)}]
    \draw [draw=dgreen, fill=lgreen] (0, 0) rectangle (4, 4);
    \node at (2,4.5) {$N_s\times N_0$};
        \node at (2,-1) {$W_p$};
    \end{scope}
    \begin{scope}[shift={(5,0)}]
        \draw [draw=dgreen, fill=lgreen] (0, 0) rectangle (1, 4);
        \node at (0.5,4.5) {$N_0\times r$};
        \node at (0.5,-1) {$V_p^\top$};
    \end{scope}
    \begin{scope}[shift={(7,0)}]
        \draw [draw=dgreen, fill=lgreen] (0, 3) rectangle (1, 4);
        \node at (0.5,4.5) {$r\times r$};
        \node at (0.5,-1) {$S_p$};
    \end{scope}
    
    \begin{scope}[shift={(9,0)}]
        \draw [draw=dgreen, fill=lgreen] (0, 3) rectangle (4, 4);
        \node at (2,4.5) {$r\times N_0$};
        \node at (5.,2) {\huge $=$};
        \node at (2,-1) {$V_p$};
    \end{scope}
    \begin{scope}[shift={(15,0)}]
        \draw [draw=dgreen, fill=lgreen] (0, 0) rectangle (1, 4);
        \node at (0.5,4.5) {$N_s\times r$};
        \node at (0.5,-1) {$U_p$};
    \end{scope}
    \begin{scope}[shift={(17,0)}]
        \draw[dgreen] (0, 3) rectangle (1, 4);
        \draw[dgreen] (0, 4) -- (1, 3);
        \node at (0.5,4.5) {$r\times r$};
        \node at (0.5,-1) {$\Sigma_p$};
    \end{scope}
    \begin{scope}[shift={(19,0)}]
        \draw [draw=dgreen, fill=lgreen] (0, 3) rectangle (4, 4);
        \node at (2,4.5) {$r\times N_0$};
        \node at (2,-1) {$V_p$};
    \end{scope}
\end{tikzpicture}
\end{center}

Since $V_p$ and $U_p$ are unitary matrices (i.e., $V_p V_p^{\hermconj}=I_r$ and $U_p^{\hermconj} U_p=I_r$), multiplying the system~(\ref{eq:vsv1}) by the matrix $V_p^{\hermconj}$ from the right and with the matrix $U_p^{\hermconj}$ from the left leads to

\begin{equation}
    \left (U_p^{\hermconj} W_p V_p^{\top} \right) S_p =  \Sigma_p,
    \label{eq:vsv2}
\end{equation}
where all three matrices $ \left (U_p^{\hermconj} W_p V_p^{\top} \right), S_p, \Sigma_p \in \R{r\times r}$. 

\begin{center}
\begin{tikzpicture}[scale=0.5]
    \begin{scope}[shift={(0,0)}]
    \draw [draw=dgreen, fill=lgreen] (0, 3) rectangle (1, 4);
    \node at (0.5,4.5) {$r\times r$};
        \node at (0.5,2) {$U_p^{\hermconj}W_pV_p^{\top}$};
    \end{scope}
    \begin{scope}[shift={(3,0)}]
        \draw [draw=dgreen, fill=lgreen] (0, 3) rectangle (1, 4);
        \node at (0.5,4.5) {$r\times r$};
        \node at (0.5,2) {$S_p$};
        \node at (2,3.5) {\huge $=$};
    \end{scope}
    \begin{scope}[shift={(6,0)}]
        \draw[dgreen] (0, 3) rectangle (1, 4);
        \draw[dgreen] (0, 4) -- (1, 3);
        \node at (0.5,4.5) {$r\times r$};
        \node at (0.5,2) {$\Sigma_p$};
    \end{scope}

\end{tikzpicture}
\end{center}

We can now solve this system for $S_p$ and ultimately reconstruct our solution of interest as $X_p = V_p^{\top}S_pV_p$ and we refer to this technique as the LR method. Similar to the previous section, we have here made the assumption that the operator, right-hand side, and unknowns share the bases $V_p$ and $U_p$. For the numerical problems in Sections~\ref{subsec:ds1} and~\ref{subsec:ds2}, this assumption is validated empirically. A theoretical analysis of the limitations of such an assumption is the subject of future research.

To summarize, we observe that the parametrization $X_p = V_p^{\top}S_p V_p$ can dramatically reduce the size of the solution (i.e., from inverting a matrix of size $N_0\times N_0$ to inverting a matrix of size $r\times r$) and, thus, the computational complexity of the associated inverse problem. Impressively, this reduction in complexity does not incur a significant compromise in the solution's quality.  Furthermore, this approach empowers us to benefit from the symmetry of the solution $X_p$, without requiring the introduction of any reciprocity preconditioner in the solution of the inverse problem~\cite{varg-mdd-2021}.

\subsection{\HT approximation}
Unfortunately, the matrices $W_p$, $X_p$, $B_p$ are usually not globally low-rank in 3-dimensional problems. However, on the bright side, those matrices can be well approximated using block low-rank structures. Particularly, they can be represented in a TLR structure~\cite{hong-tlr-2021} as well as \HT structure~\cite{hack-h2-2000,grro-fmm-1987}, as we verify it in our numerical experiments in Sections~\ref{sec:ds3}.

In this section, we introduce the fundamental concepts of the \HT matrix approximation. For more comprehensive and formal definitions, please refer to~\cite{grro-fmm-1987,hack-h2-2000,borm-h2-2010}. It is important to note that our version for \HT decomposition differs from the standard one, yet it remains valid. A thorough explanation of our version of \HT is available in~\cite{sush-sne-2020}.

Let us, for example, consider the dense right-hand side matrix $B_p$ in the system~(\ref{eq:fr_sys}) and build \HT factorization of it. Noting that rows and columns of $B_p$ are associated with sources and receivers in the physical acquisition geometry, we can start by hierarchically splitting the grids of sources and receivers and storing these divisions within tree structures. Utilizing this grid division, we rearrange the rows and columns of the matrix $B_p$: such a rearrangement creates a block structure in the new matrix, which possesses now block low-rank properties. Let the matrix $B_p$ have $N_{r}$ block rows and $N_{c}$ block columns. Thus, the matrix $B_p$ is divided into blocks 
$$B_{ij}\in \R{N_{i} \times N_{j}}, \quad i \in 1\dots N_{r}, j \in 1\dots N_{c},$$ 
where $N_{i}$ is a size of $i$-th block row, and $N_{j}$ is a size of $j$-th block column.
We also define block rows and block columns
$$B_{i}\in \R{N_{i} \times N_{0}},$$
$$B_{j}\in \R{N_{s} \times N_{j}}.$$

\begin{figure}[H]
\begin{subfigure}{0.9\textwidth}
\begin{center}
\begin{tikzpicture}[scale=0.5]
    \begin{scope}[shift={(0,0)}]
        \fill[lgreen] (0, 0) rectangle (6, 1);
        \fill[dgreen] (0, 0) rectangle (1, 1);
        \fill[dgreen] (2, 0) rectangle (3, 1);
        \draw[step=1cm,dgreen] (0,0) grid (6,1);
    \end{scope}
\end{tikzpicture}
\end{center}
\caption{Illustration of $B_i$} 
\end{subfigure}\hspace*{\fill}

\medskip
\begin{subfigure}{0.9\textwidth}
\begin{center}
\begin{tikzpicture}[scale=0.5]
    \begin{scope}[shift={(0,0)}]
        \fill[lgreen] (0, 0) rectangle (6, 1);
        \fill[white] (0, 0) rectangle (1, 1);
        \fill[white] (2, 0) rectangle (3, 1);
        \draw[step=1cm,dgreen] (0,0) grid (6,1);
    \end{scope}
\end{tikzpicture}
\end{center}
\caption{Illustration of $\widehat{B}_i$} 
\end{subfigure}\hspace*{\fill}
\caption{Example of $B_i$ and $\widehat{B}_i$. Dark green blocks represent close blocks and light green blocks represent far blocks. White blocks are zeros.} 
\label{fig:h2_bi}
\end{figure}
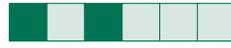
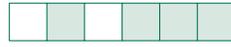

All blocks in matrix $B_p$ can be separated into two classes: blocks that correspond to the sources and receivers that are located closer than some predefined distance (called 'close blocks`) and all others (called 'far blocks`). Let us also define the block row $\widehat{B}_{i}$, which is $B_i$ with zeroed close blocks, and accordingly, block column $\widehat{B}_{j}$ is $B_j$ with zeroed close blocks:

$$ B_{ij} \in \widehat{B}_i =  \begin{cases}
    B_{ij} \in B_i       & \quad \mathbf{if~}  B_{ij} \mathbf{~far} \\
    0                & \quad \mathbf{if~} B_{ij} \mathbf{~close} 
    \end{cases}$$ 

$$ B_{ij} \in \widehat{B}_j =  \begin{cases}
B_{ij} \in B_j       & \quad \mathbf{if~}  B_{ij} \mathbf{~far} \\
0                & \quad \mathbf{if~} B_{ij} \mathbf{~close} 
\end{cases}$$
See illustration of $B_i$ and $\widehat{B}_i$ in
Figure~\ref{fig:h2_bi}.

\begin{figure}[H]
\begin{center}
\begin{tikzpicture}[scale=0.5]
    \begin{scope}[shift={(0,0)}]
        \fill[lgreen] (0, 0.5) rectangle (6, 1);
        \fill[dgreen] (0, 0) rectangle (1, 1);
        \fill[dgreen] (2, 0) rectangle (3, 1);
        \draw[step=1cm,dgreen] (0,0) grid (6,1);
    \end{scope}
\end{tikzpicture}
\end{center}
\caption{Example of block row $U_i^{\hermconj} B_i$.} 
\label{fig:h2_bi_t}
\end{figure}
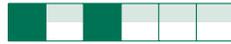

The key property of the \HT matrix formulation used in this work is that block row $\widehat{B}_{i}$ and block column $\widehat{B}_{j}$ matrices are globally low rank. In other words, if we consider the left factor $U_i\in\R{N_i\times N_i}$ of the SVD decomposition of the $\widehat{B}_{i}$ matrix, and multiply the block row $\widehat{B}_{i}$ from the left with $U_i^{\hermconj}$, the resulting matrix will contain non-zero values only in the first $r$ rows (whilst all other rows will have values close to zero), where $r$ is rank of $\widehat{B}_{i}$. Similarly, if we multiply the block row $B_i$ from the left with $U_i^{\hermconj}$, the resulting matrix will have the same compression for the far blocks and some changes in the close blocks.  
See the illustration of the block row $U_i^{\hermconj} B_i$ in Figure~\ref{fig:h2_bi_t}.

Analogously, we can compress the block columns by computing the right SVD factor $V_j\in\R{N_j\times N_j}$ of the block columns $\widehat{B}_j$ and multiplying  $V_j$ from the right to the block column $B_j$. 

Combining all the obtained $U_i$ and $V_j$ blocks into two big block-diagonal unitary matrices, obtain:
 $$U = \begin{bmatrix}
    U_{1}&&\\
    &\ddots&\\
    &&U_{N_b}\\
\end{bmatrix}, \quad
V = \begin{bmatrix}
    V_{1}&&\\
    &\ddots&\\
    &&V_{N_b}\\
\end{bmatrix},$$
where $U\in \R{N_s \times N_s}$, $V\in \R{N_0 \times N_0}.$ Note, that in practice, it is more efficient to obtain $V_j$ factor computing SVD not of the original $\widehat{B}_{j}$ block columns, but of the block column $U^{\hermconj}B_j$ with excluded zeros. 

Multiplying the the original matrix $B_p$ from left and right with the matrices $U^{\hermconj}$ and $V^{\hermconj}$ we obtain a sparse matrix $S_B\in\R{N_s\times N_0}$:
\begin{equation}
    \label{eq:subv}
   S_B = U^{\hermconj} B_p V^{\hermconj}. 
\end{equation}
More specifically, $S_B$ is a block matrix with the following block structure: it has a dense block $S_{ij} = U_{i}^{\hermconj}B_{ij}V^{\hermconj}_{i}$, if the corresponding block of $B_p$ is of close type and it has zero block with non-zero $r\times r$ upper-left corner, if the corresponding block of $B_p$ is of far type. Note that in the actual implementation of our algorithm, we store the matrix $S_B$ in a block-sparse format to avoid storage of the zero elements.
Multiplying the Equation~(\ref{eq:subv}) by the matrices $U$ and $V$, obtain:
$$B_p = US_BV$$
See Figure~\ref{fig:h2} for an illustration of the structure of our matrix's parametrization:

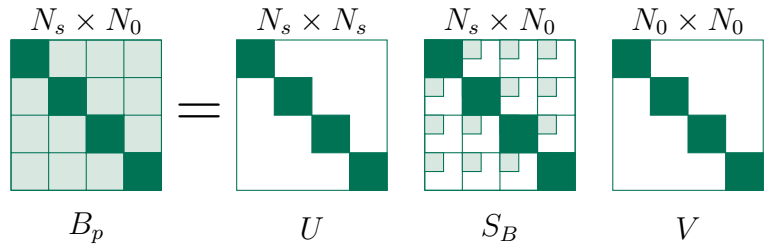
\begin{figure}[H]
\begin{center}
\begin{tikzpicture}[scale=0.5]
    \begin{scope}[shift={(0,0)}]
    \draw [draw=dgreen, fill=lgreen] (0, 0) rectangle (4, 4);
    \node at (2,4.5) {$N_s\times N_0$};
        \node at (2,-1) {$B_p$};
        \draw[step=1.0,dgreen,thin] (0,0) grid (4,4);
        \foreach \x in {0,...,3}{
            \draw [draw=dgreen,fill=dgreen] (\x, 3-\x) rectangle (\x+1, 4-\x);
        }
        \node at (5.,2) {\huge $=$};
    \end{scope}
    \begin{scope}[shift={(6,0)}]
        \draw [draw=dgreen] (0, 0) rectangle (4, 4);
        \foreach \x in {0,...,3}{
            \draw [draw=dgreen,fill=dgreen] (\x, 3-\x) rectangle (\x+1, 4-\x);
        }
        \node at (2,4.5) {$N_s\times N_s$};
        \node at (2,-1) {$U$};
    \end{scope}
    \begin{scope}[shift={(11,0)}]
        \draw [draw=dgreen] (0, 0) rectangle (4, 4);
        \draw[step=1.0,dgreen,thin] (0,0) grid (4,4);
        \foreach \x in {0,...,3}
        \foreach \y in {3,...,0}{
            \draw [draw=dgreen,fill=lgreen] (\x, \y+0.5) rectangle (\x+0.5, \y+1);
        }
        \foreach \x in {0,...,3}{
            \draw [draw=dgreen,fill=dgreen] (\x, 3-\x) rectangle (\x+1, 4-\x);
        }
        \node at (2,4.5) {$N_s\times N_0$};
        \node at (2,-1) {$S_B$};
    \end{scope}
    \begin{scope}[shift={(16,0)}]
        \draw [draw=dgreen] (0, 0) rectangle (4, 4);
        \foreach \x in {0,...,3}{
            \draw [draw=dgreen,fill=dgreen] (\x, 3-\x) rectangle (\x+1, 4-\x);
        }
        \node at (2,4.5) {$N_0\times N_0$};
        \node at (2,-1) {$V$};
    \end{scope}
\end{tikzpicture}
\end{center}
\caption{Illustration of \HT structure of matrix $B_p$.}
\label{fig:h2}
\end{figure}
Note that this is the one-level version of the \HT matrix. The multi-level version is more advanced, and its application is the subject of future research.

\subsection{$\mathcal{H}^2$-like parametrization of the MDD}
In this section, to simplify our notation, we drop the $p$ index from every matrix. From now on, we therefore assume that the equations hold true for all frequencies $p\in\R{N_f}$.

As in Section~\ref{sec:xvsv}, we know {\em a priori} that $X$ is symmetric (but not hermitian). We also want the matrices $B$ and $X$ to share the right factor. Thus, we will search $X$ in the following format: 
$$X = V^{\top}S_X V,$$
where $S_X$ has the same format as $S_B$. We obtain the system:
$$ WV^{\top}S_X V = US_B V.$$
Then we apply the matrices $U\hermconj$ to the $W$ from the left and $V\hermconj$ from the right. We obtain the $U\hermconj WV\hermconj$ operator. Thanks to the physical properties of the problem, that operator has the same structure as $S_B$:
$$S_W = U\hermconj WV\hermconj.$$
We obtain the following system:
\begin{equation} \label{eq:sys1}
US_W\conj{V} V^{\top} S_X V = U S_B V.
\end{equation}
Let us then multiply the system (\ref{eq:sys1}) with $U\hermconj$ from the left and $V\hermconj$ from the right:
$$U\hermconj US_W\conj{V} V^{\top} S_X V V\hermconj= U\hermconj U S_B VV\hermconj.$$
Considering the fact that $U$ and $V$ are hermitian matrices, and thus $U\hermconj U = I$,  $VV\hermconj = I$, $\conj{V} V^{\top} = (VV\hermconj)^{\top} = I^{\top}=I$, we receive:
\begin{equation}
\label{eq:sss}
S_W S_X = S_B.    
\end{equation}
We can now solve this system with the iterative method LSQR~\cite{paig-lsqr-1982} for $S_X$.
After computing $S_X$, we reconstruct the solution using the formula
$$X =  V^{\top} S_X V.$$
After that, using the Fourier transform, we obtain the final solution in the time domain. 

\begin{remark}
While solving the system~(\ref{eq:sss}), the LSQR method employs matrix-by-matrix multiplication. If the process does not preserve sparsity, the resulting dense matrix negates the benefits of compression. The solution is to put to zero all the elements of the resulting matrix that were zero in the initial matrices.
This approach is illustrated in Figure~\ref{fig:fill}, which visually explains how we preserve the sparse structure during multiplication.
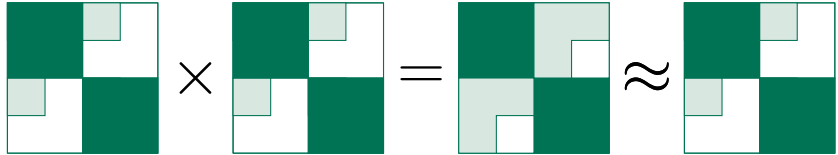
\begin{figure}[H]
\begin{center}
\begin{tikzpicture}[scale=0.5]
    \begin{scope}[shift={(-1,0)}]
    \node at (5.,2) {\huge $\times$};
    \draw [draw=dgreen] (0, 0) rectangle (4, 4);
    \draw[step=2.0, dgreen,thin] (0,0) grid (4,4);
    \foreach \x in {0,...,1}
        \foreach \y in {0,...,1}{
            \draw [draw=dgreen,fill=lgreen] (2*\x, 2*\y+1) rectangle (2*\x+1, 2*\y+2);
        }
    \foreach \x in {0,...,1}{
            \draw [draw=dgreen,fill=dgreen] (2*\x, 2-\x*2) rectangle (2*\x+2, 4-\x*2);
        }
    \end{scope}
    \begin{scope}[shift={(5,0)}]
        \draw [draw=dgreen] (0, 0) rectangle (4, 4);
        \draw[step=2.0, dgreen,thin] (0,0) grid (4,4);
        \foreach \x in {0,...,1}
        \foreach \y in {0,...,1}{
                \draw [draw=dgreen,fill=lgreen] (2*\x, 2*\y+1) rectangle (2*\x+1, 2*\y+2);
        }
        \foreach \x in {0,...,1}{
                \draw [draw=dgreen,fill=dgreen] (2*\x, 2-\x*2) rectangle (2*\x+2, 4-\x*2);
            }
        \node at (5.,2) {\huge $=$};
    \end{scope}
    \begin{scope}[shift={(11,0)}]
        \draw [draw=dgreen] (0, 0) rectangle (4, 4);
        \draw[step=2.0, dgreen,thin] (0,0) grid (4,4);
    \foreach \x in {0,...,1}{
            \draw [draw=dgreen,fill=dgreen] (2*\x, 2-\x*2) rectangle (2*\x+2, 4-\x*2);
    }
    \foreach \x in {0,...,1}{
            \draw [draw=dgreen,fill=lgreen] (2*\x, \x*2) rectangle (2*\x+2, \x*2+2);
    }
    \foreach \x in {0,...,1}{
            \draw [draw=dgreen,fill=white] (2*\x+1, \x*2) rectangle (2*\x+2, \x*2+1);
    \node at (5.,2) {\huge $\approx$};
    }
    \end{scope}
    \begin{scope}[shift={(17,0)}]
    \draw [draw=dgreen] (0, 0) rectangle (4, 4);
    \draw[step=2.0, dgreen,thin] (0,0) grid (4,4);
    \foreach \x in {0,...,1}
        \foreach \y in {0,...,1}{
            \draw [draw=dgreen,fill=lgreen] (2*\x, 2*\y+1) rectangle (2*\x+1, 2*\y+2);
        }
    \foreach \x in {0,...,1}{
            \draw [draw=dgreen,fill=dgreen] (2*\x, 2-\x*2) rectangle (2*\x+2, 4-\x*2);
        }
        
    \end{scope}
\end{tikzpicture}
\end{center}
\caption{The fill-in during the matrix multiplication.}
\label{fig:fill}
\end{figure}
\end{remark}

\begin{remark}
In the numerical implementation, matrix-by-matrix multiplication is carried out in a specific way: the right matrix is converted into a vector that contains only its nonzero elements. The left matrix is implemented as a linear operator, which is then applied to the vectorized right matrix. After applying the linear operator, the resulting vector is converted back into a sparse matrix format. This approach efficiently handles the multiplication by using only the nonzero elements. Another significant benefit of this method is its compatibility with the standard LSQR algorithm, which naturally operates with vectors and linear operators. This compatibility makes integrating the LSQR method with this approach easier.
\end{remark}

 This approach empowers us to benefit from the symmetry of the solution $X$, enabling, for example, the utilization of the reciprocity precondition. 
 We refer to this technique as the \HT method.

\section{Numerical Examples}
In this section, we present three numerical experiments to validate the proposed methods; these involve applying global low-rank (for the 2-dimensional case) and block low-rank (for the 3-dimensional case) compression techniques to all elements of the linear system we want to solve. The different experiments are chosen to mirror real-world seismic exploration scenarios. 
Among the datasets, two are two-dimensional, offering standard complexity scenarios, while the third is a significantly larger three-dimensional dataset, providing a more challenging and comprehensive test of the capabilities of our approaches in handling large-scale seismic data. All algorithms are implemented in Python; our focus is to characterize the time and memory savings of our methods in return for a controllable compromise in the quality of the results.

\subsection{Ocean-bottom cable (2d)}
\label{subsec:ds1}
In this subsection, we consider a synthetic ocean-bottom cable (OBC) dataset. In our specific case study, we use a two-dimensional segment of the SEG/EAGE Overthrust model~\cite{amin-over-1997}, which includes an additional water layer on top. 401 receivers are arranged on the sea floor. Each is spaced 10 meters apart and 292 meters below the ocean's surface. 401 sources are placed at a depth of 20 m with a lateral spacing of 10 m.

\begin{figure}[H]
\begin{subfigure}{0.34\textwidth}
\includegraphics[width=\linewidth]{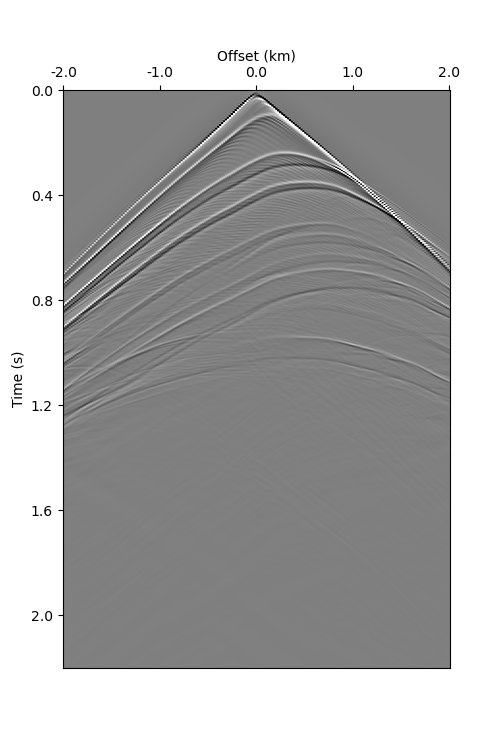}
\caption{Modeled solution} \label{fig:ds1_an}
\end{subfigure}\hspace*{\fill}
\begin{subfigure}{0.34\textwidth}
\includegraphics[width=\linewidth]{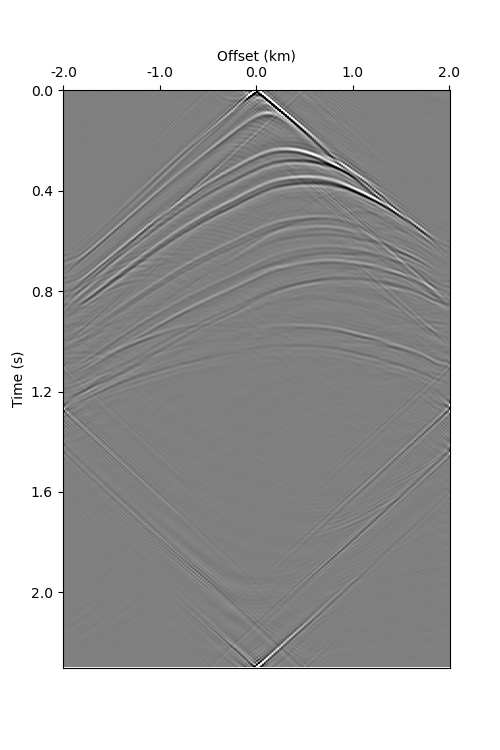}
\caption{Full, time dom.} \label{fig:ds1_full}
\end{subfigure}\hspace*{\fill}
\begin{subfigure}{0.34\textwidth}
\includegraphics[width=\linewidth]{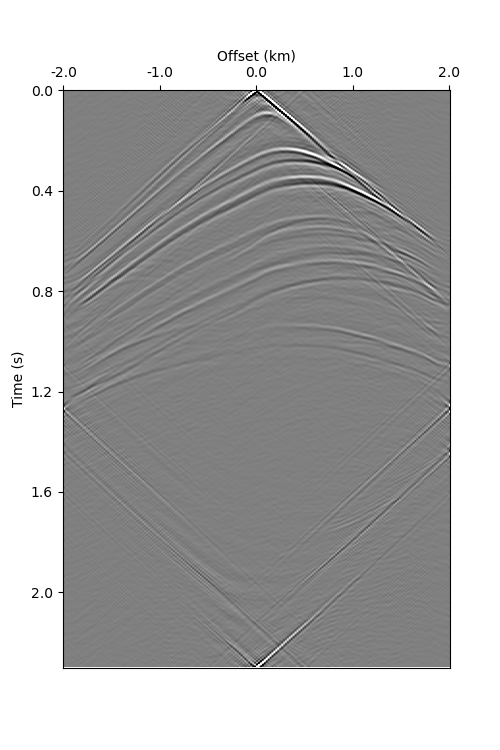}
\caption{LR method} \label{fig:ds1_lr}
\end{subfigure}
\medskip
\begin{subfigure}{0.34\textwidth}
\includegraphics[width=\linewidth]{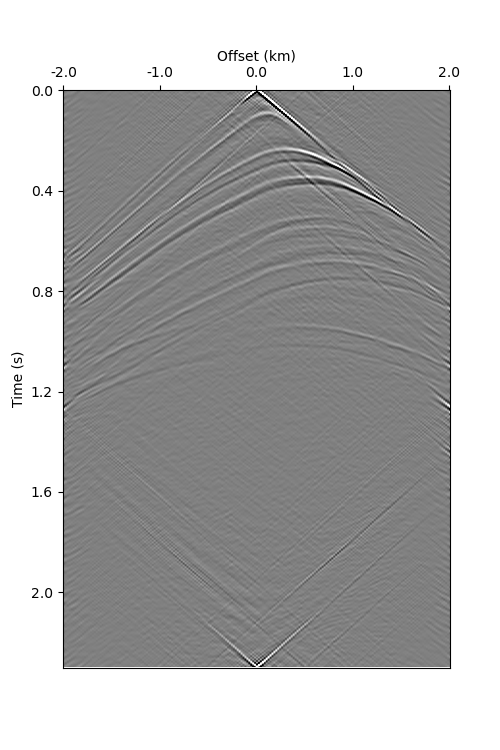}
\caption{USV method} \label{fig:ds1_usv}
\end{subfigure}\hspace*{\fill}
\begin{subfigure}{0.34\textwidth}
\includegraphics[width=\linewidth]{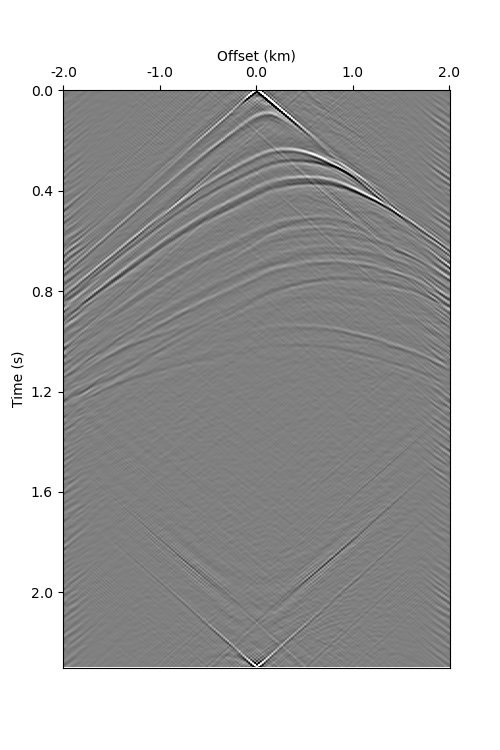}
\caption{USV-rsp method} \label{fig:ds1_suv_rcp}
\end{subfigure}\hspace*{\fill}
\begin{subfigure}{0.34\textwidth}
\includegraphics[width=\linewidth]{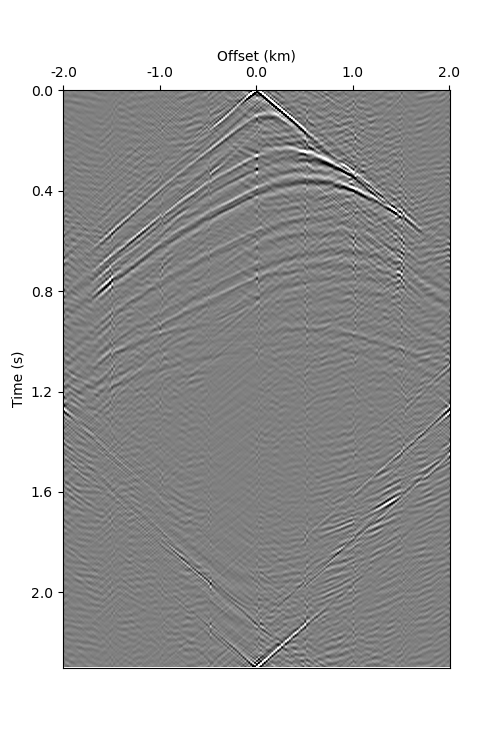}
\caption{\HT method} \label{fig:ds1_h2}
\end{subfigure}
\caption{MDD inversion results for different solution methods.}
\label{fig:ds1}
\end{figure}

Figure~\ref{fig:ds1} illustrates the inversion results achieved using MDD, using the LR, USV, USV-rsp (USV with reciprocity precondition), and \HT compression strategies described in the paper. These results are further compared to the reference response directly modeled via finite difference without the water column and free surface and to the "full" solution where we solve MDD in the time or frequency domain with uncompressed matrices. 

The LR method outperforms the USV and USV-rsp methods, using a significantly larger amount of data in the solution process. In this context, applying reciprocity preconditioning has minimal impact on the outcome. In the 2D scenario, the \HT method appears to be the least effective. Its complex data structure necessitates large computational time, and the additional approximations result in a visible accuracy loss. However, in 3D, the large size of the matrices involved makes this format advantageous, offsetting the initial drawbacks with its benefits, as we will show later.

Table~\ref{tab:ds1} compares the time and memory requirements for the MDD methods described above. We consider the uncompressed time and frequency domain MDD, LR, USV, and \HT methods. For timing, we present integral solution timings, including the approximation to the LR, USV, or \HT format. For memory, we sum up the sizes of all matrices involved in the MDD process (operator, right-hand side, and unknown matrices) in the frequency domain. 

Below, we present the formulas for the memory requirements. For this 2d ocean-bottom cable dataset, the sizes of the arrays are: 
number of sources $N_s = 201$, number of receivers and virtual receivers $N_0 = 201$, number of time samples $N_t = 1150$, and  number of frequency samples $N_f = 200$, rank is set to $r=50$.

The full-matrix right-hand side data in the time domain has a size of
\begin{equation}
    \label{eq:mem_rhs_t}
    N_s\times N_0 \times N_t \times 4 \text{Bytes} = 185.8\text{MB},
\end{equation}
whilst the size of the full-matrix right-hand side data in the frequency domain is
\begin{equation}
    \label{eq:mem_rhs_f}
    N_s\times N_0\times N_f \times 8 \text{Bytes} = 64.6\text{MB},
\end{equation}
the size of the LR right-hand side data is
\begin{equation}
    \label{eq:mem_rhs_lr}
    N_s\times r \times N_f \times 8 \text{Bytes} = 16.1\text{MB},
\end{equation}
the size of the USV right-hand side data is
\begin{equation}
    \label{eq:mem_rhs_usv}
    r \times r \times N_f \times 8 \text{Bytes} = 4.0\text{MB},
\end{equation}
and the size of the \HT right-hand side data is computed empirically as $13.5\text{MB}$.

The full-matrix unknowns in the time domain have a size of
\begin{equation}
    \label{eq:mem_x_t}
    N_0\times N_0\times N_t \times 4 \text{Bytes} = 185.8\text{MB},
\end{equation}
whilst the size of the full-matrix unknowns in the frequency domain is 
\begin{equation}
    \label{eq:mem_x_f}
    N_0\times N_0\times N_f \times 8 \text{Bytes} = 64.6\text{MB},
\end{equation}
the size of the LR unknowns is 
\begin{equation}
    \label{eq:mem_x_lr}
    N_0\times r\times N_f \times 8 \text{Bytes} = 16.1\text{MB},
\end{equation}
the size of the USV unknowns is 
\begin{equation}
    \label{eq:mem_x_usv}
    r\times r\times N_f \times 8 \text{Bytes} = 4.0\text{MB},
\end{equation}
and the size of the \HT unknowns side data is again computed empirically as $13.5\text{MB}$.

The operator is always stored in the frequency domain. The operators have the same size as the frequency domain right-hand side of the same method for all methods except for LR, where operators have the same size as the full-matrix right-hand side.

\begin{table}[H]
\centering
\begin{tabular}{c||c|c|c|c|cc}
 &  full, time dom. & full, freq. dom & LR & USV & \HT \\ \hline \hline
Time (s) & 21.3 & 12.8 & 6.0 & 0.71 & 33.8\\
Mem (MB)            &       436.3      & 193.9  &   96.8  &   12.0  &  40.4\\
\end{tabular}
\caption{Time and memory comparison of the candidate methods.}
\label{tab:ds1}
\end{table}

In this dataset, all three system matrices across every frequency possess a globally low rank. This characteristic underscores the effectiveness of the USV method, which emerges as the optimal choice in terms of both memory and time efficiency. 

\subsection{Target level below a complex salt body (2d)}
\label{subsec:ds2}
\begin{figure}[H]
\includegraphics[width=12cm]{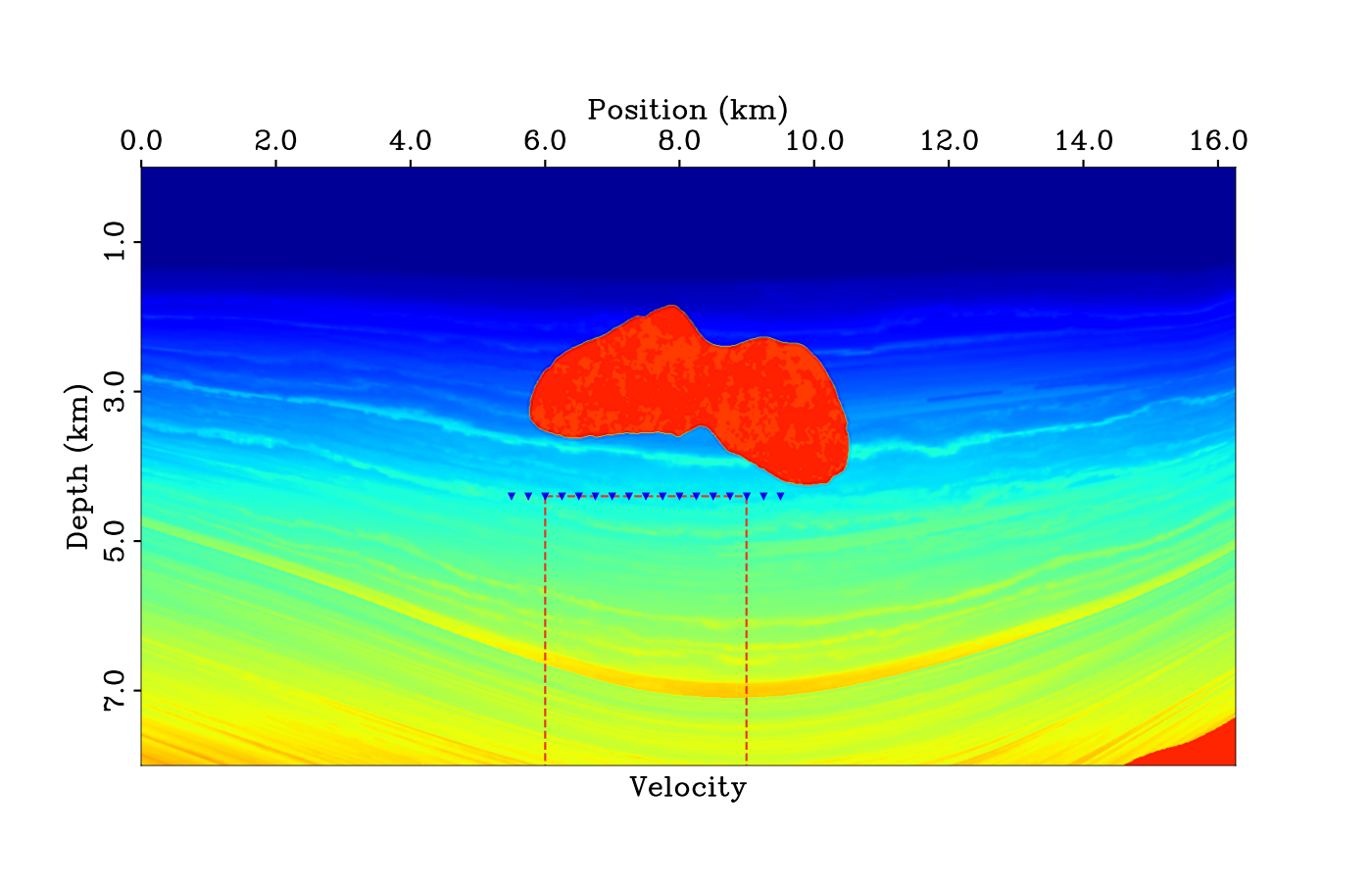}
\caption{Sub-salt velocity model supporting decomposed wavefields at the receiver array. A magenta line runs horizontally at a depth of 4.4 kilometers below the layer above the salt. The sources of the wavefields, shown as a solid red line, are positioned at the surface. A dashed black line also marks the specific target area where we aim to capture and analyze the reflected wave signals.}
\label{fig:salt}
\end{figure}

For our second experiment, we use a dataset created by redatuming surface reflection data at a target level below a complex salt body~\cite{varg-mdd-2021}. These wavefields are obtained by applying the scattering Rayleigh-Marchenko method as detailed in~\cite{varg-salt-2021}. This method allows one to estimate wavefields traveling upwards and downwards as if they were emitted from surface sources and then captured by virtual receivers located under a complex layer of earth, like beneath a salt body; see Figure~\ref{fig:salt}.

The experimental area extends over 16.26 km horizontally and reaches a depth of 8 km. It is designed to imitate a challenging environment with an uneven salt layer and other complex features like sharp changes in rock layers, uneven sediment layers, and inconsistencies between different rock layers. Redatuming in such a complex area is difficult because of the stark differences in rock properties, which can cause multiple internal reflections. These reflections can lead to misleading signals in the data, often seen as artifacts or distortions. Additionally, other challenges like limited frequency range and insufficient source illumination make the redatuming process more challenging.

\begin{figure}[H]
\begin{subfigure}{0.34\textwidth}
\includegraphics[width=\linewidth]{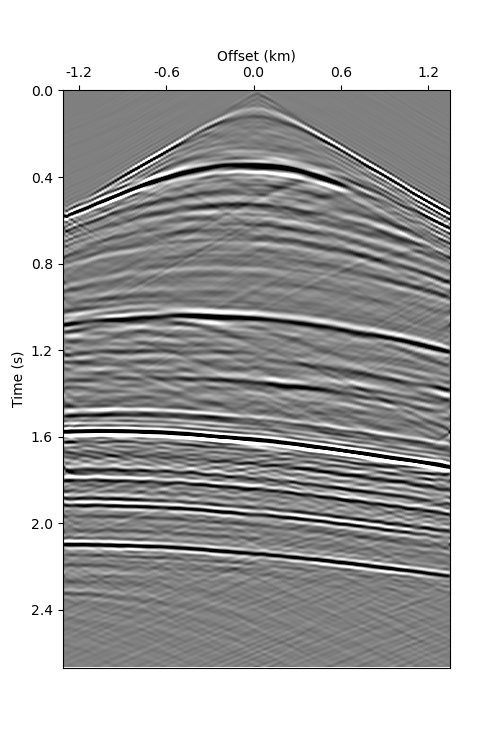}
\caption{Modeled solution} \label{fig:ds2_an}
\end{subfigure}\hspace*{\fill}
\begin{subfigure}{0.34\textwidth}
\includegraphics[width=\linewidth]{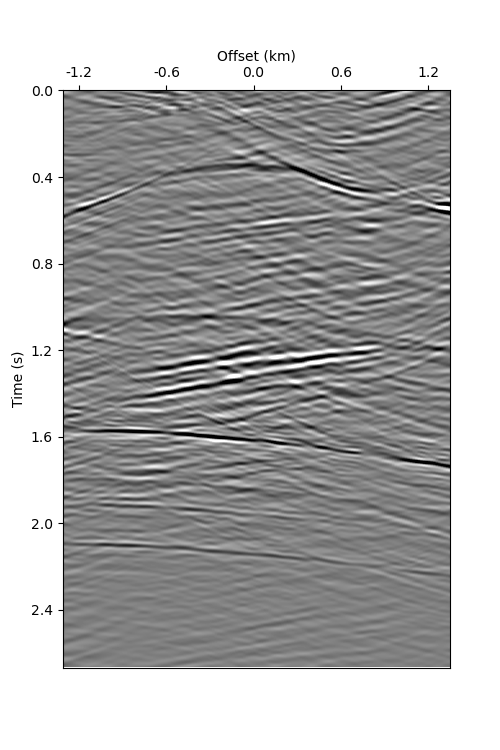}
\caption{Full, time dom} \label{fig:ds2_full}
\end{subfigure}\hspace*{\fill}
\begin{subfigure}{0.34\textwidth}
\includegraphics[width=\linewidth]{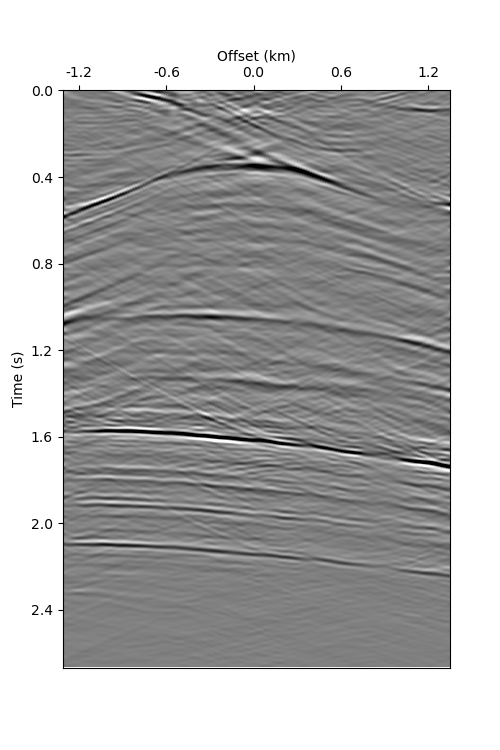}
\caption{Full-rcp, time dom} \label{fig:ds2_full_s}
\end{subfigure}
\caption{MDD inversion results for reference solution.} 
\label{fig:ds2_ref}
\end{figure}

\begin{figure}[H] 
\begin{subfigure}{0.34\textwidth}
\includegraphics[width=\linewidth]{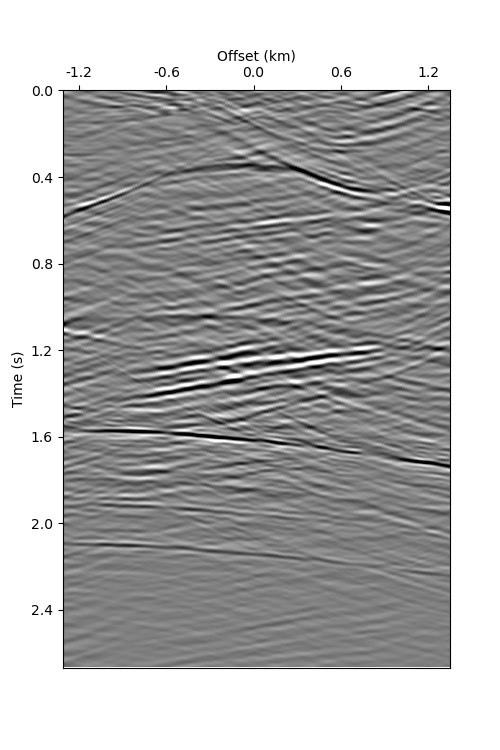}
\caption{LR method} \label{fig:ds2_lr}
\end{subfigure}\hspace*{\fill}
\begin{subfigure}{0.34\textwidth}
\includegraphics[width=\linewidth]{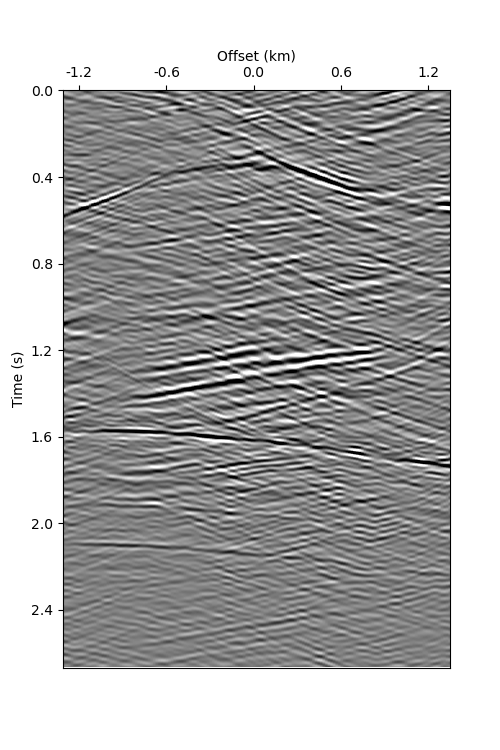}
\caption{USV method} \label{fig:ds2_usv}
\end{subfigure}\hspace*{\fill}
\begin{subfigure}{0.34\textwidth}
\includegraphics[width=\linewidth]{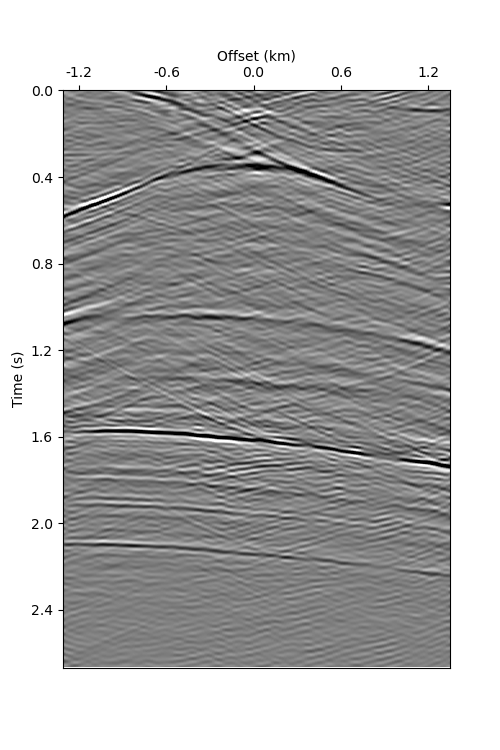}
\caption{USV-rsp method } \label{fig:ds2_usv_s}
\end{subfigure}
\medskip

\centering
\begin{subfigure}{0.34\textwidth}
\includegraphics[width=\linewidth]{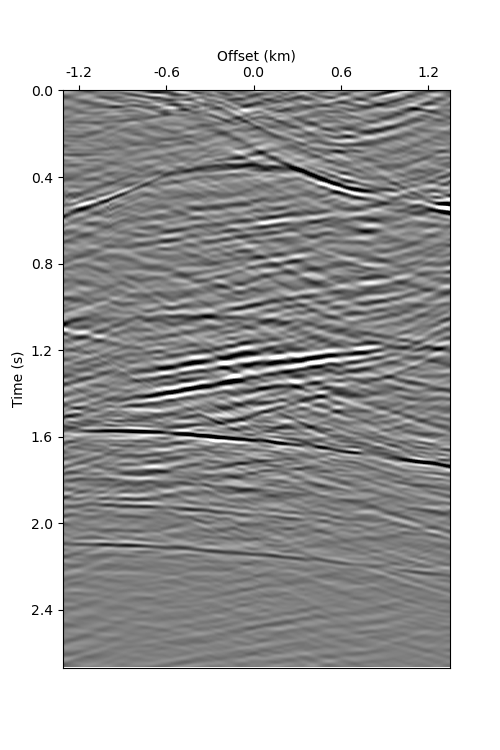}
\caption{\HT method} \label{fig:ds2_h2}
\end{subfigure}
\begin{subfigure}{0.34\textwidth}
\includegraphics[width=\linewidth]{figs/ds2_h2.png}
\caption{$\mathcal{H}^2$-rsp method} \label{fig:ds2_h2_s}
\end{subfigure}
\caption{MDD inversion results for different solution methods.} 
\label{fig:ds2_sol}
\end{figure}

Figure~\ref{fig:ds2_ref} presents a reference solution, and Figure~\ref{fig:ds2_sol} presents different results obtained by applying MDD inversion to the considered dataset. These results were achieved with the various methods detailed in the paper. This particular dataset benefits from having a globally low-rank structure across all three matrices in systems for all frequencies. However, despite this advantage, the problem is inherently unstable. The results are only reliable when the reciprocity preconditioner is applied, as clearly demonstrated in Figure~\ref{fig:ds2_sol}. (Compare same methods with and without preconditioning.)


Table~\ref{tab:ds2_mem} compares the time and memory requirements for the MDD methods described above. We consider the uncompressed time and frequency domain MDD, LR, USV, and \HT methods. For timing, we present integral solution timings, including the approximation to the LR, USV, or \HT format. For memory, we sum up the sizes of all matrices involved in the MDD process (operator, right-hand side, and unknown matrices) in the frequency domain. 

We compute the memory requirements using the formulas~(\ref{eq:mem_rhs_t} --
\ref{eq:mem_x_usv}) presented in the previous subsection. The sizes of the \HT unknowns side data are computed empirically. In this 2d target level below a complex salt body dataset, the number of sources is $N_s = 201$, the number of receivers and virtual receivers is $N_0 = 151$, the number of time samples is $N_t = 2001$, the number of frequency samples is $N_f = 400$, and rank is set to $r=50$.

\begin{table}[H]
\centering
\begin{tabular}{c||c|c|c|c|cc}
 & full, time dom. & full, freq. dom. & LR & USV & \HT \\ \hline \hline
Time (s) & 26.5 & 19.3 & 10.8 & 2.1 & 19.5 \\
Mem (MB)            &       583.0    &    267.2         &   153.4  &   24.0  &  80.7\\
\end{tabular}
\caption{Time and memory comparison of the candidate methods.}
\label{tab:ds2_mem}
\end{table}

For this dataset, the reciprocity preconditioning is crucial to obtaining valuable results. Thus, the USV-rsp method outperforms other methods not only in terms of memory and time but also in terms of solution quality.
The LR method's inability to use reciprocity becomes critical in this example. As in the previous example,  the \HT method is the least effective for the above reasons.

\subsection{Ocean-bottom cable (3d)}
\label{sec:ds3}
In this subsection, we consider the synthetic ocean-bottom cable model. In this example, we use a full 3d  SEG/EAGE Overthrust model~\cite{amin-over-1997}, which includes an additional water layer on top. 

The grid of $177\times 90$ receivers is arranged on the sea floor. Each spaced 20 meters apart and 292 meters below the ocean's surface. Additionally, we create a $120\times 217$ grid of shot gathers with a source spacing of 20 m and a depth of 10 m.

\begin{figure}[H]

\begin{subfigure}{0.9\textwidth}
\includegraphics[width=\linewidth]{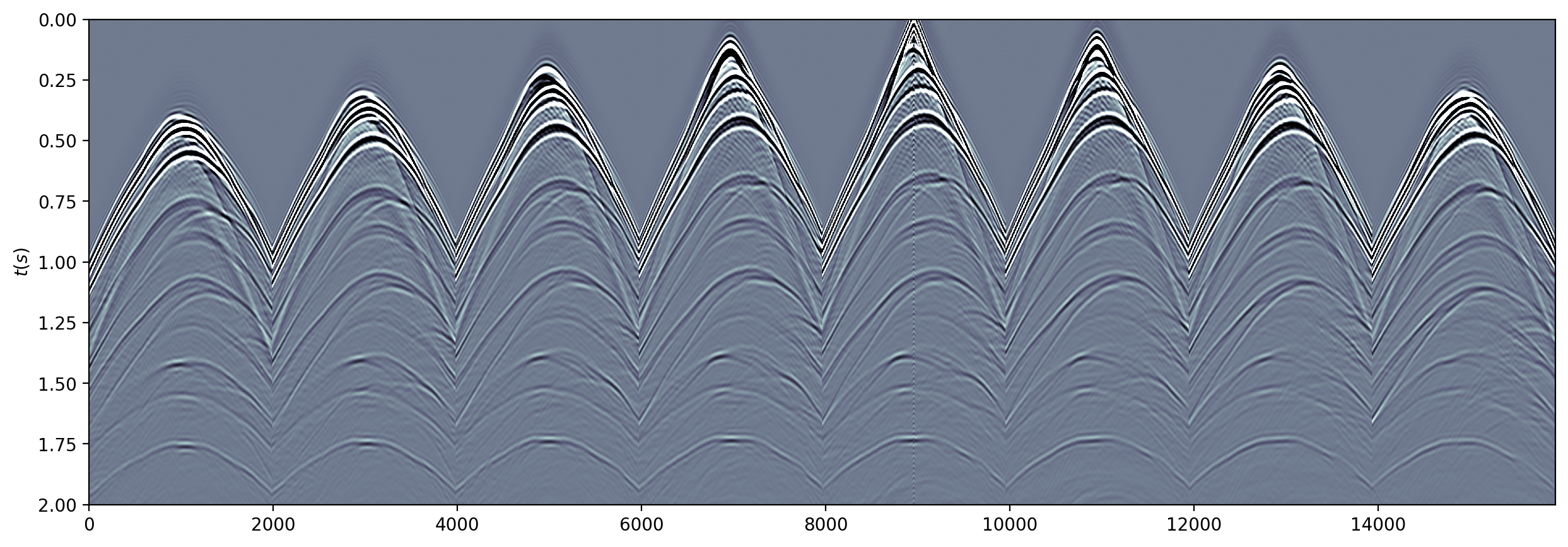}
\caption{Modeled local reflectivity} \label{fig:ds3_an}
\end{subfigure}\hspace*{\fill}

\medskip
\begin{subfigure}{0.9\textwidth}
\includegraphics[width=\linewidth]{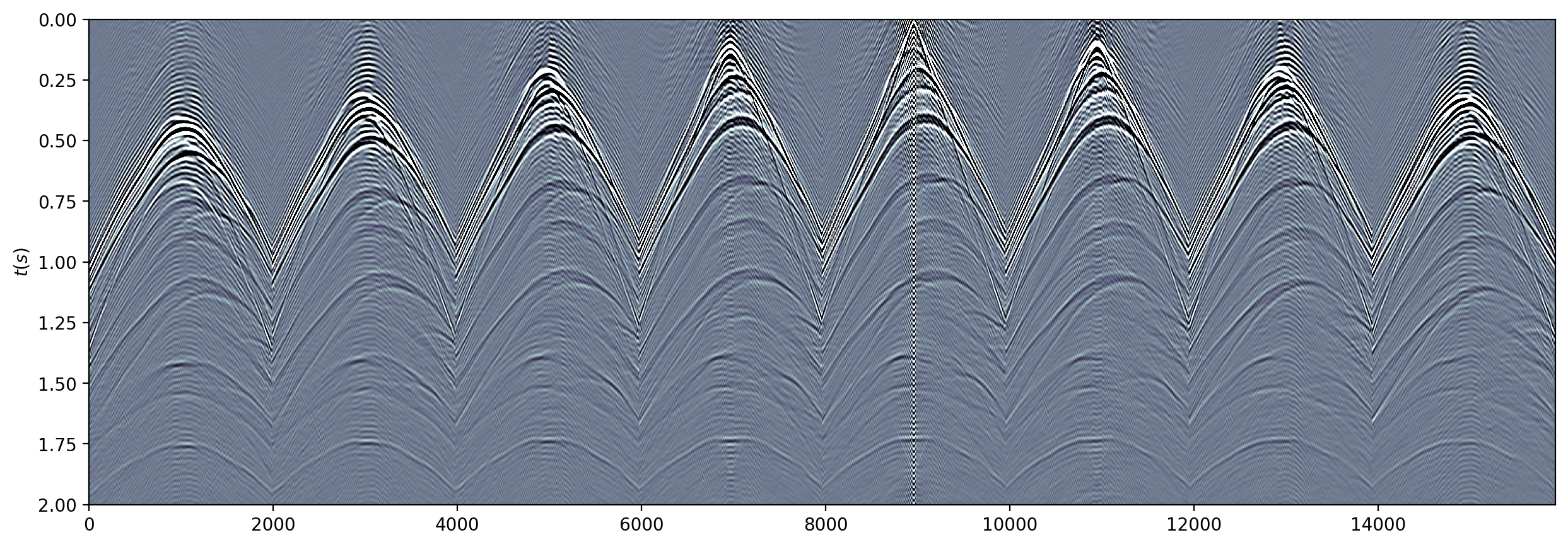}
\caption{Solution $\mathbf{R_{inv}}$, computed in frequency domain, for full matrices} \label{fig:ds3_full}
\end{subfigure}\hspace*{\fill}
\caption{Reference solutions} 
\label{fig:ds3_r}
\end{figure}

\begin{figure}[H] 
\begin{subfigure}{0.9\textwidth}
\includegraphics[width=\linewidth]{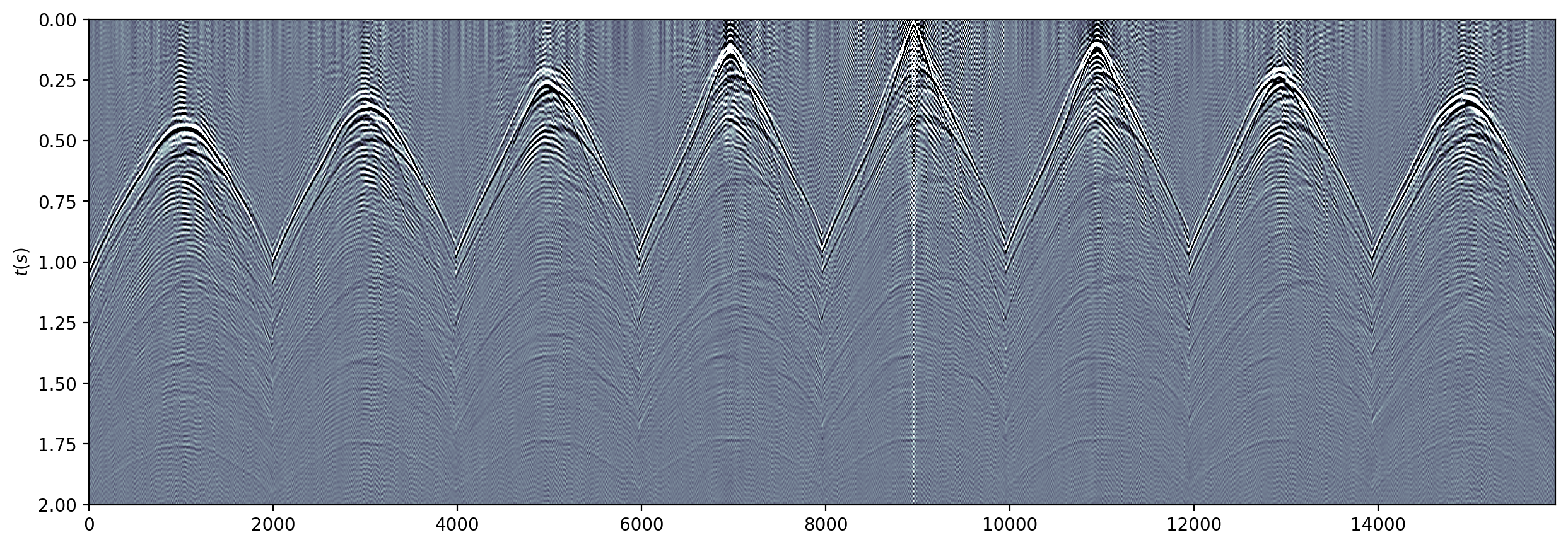}
\caption{$\mathbf{R_{inv}}$, \HT, $r=200$} \label{fig:ds3_r200}
\end{subfigure}\hspace*{\fill}

\medskip
\begin{subfigure}{0.9\textwidth}
\includegraphics[width=\linewidth]{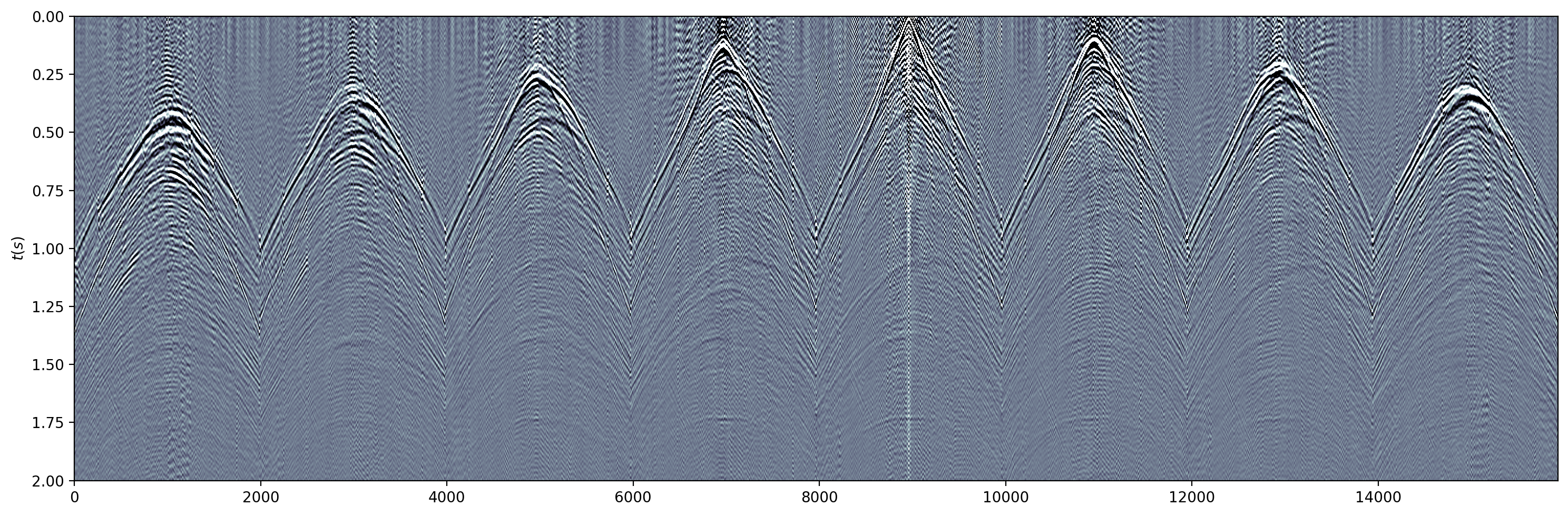}
\caption{$\mathbf{R_{inv}}$, \HT, $r=100$} \label{fig:ds3_r100}
\end{subfigure}\hspace*{\fill}

\medskip
\begin{subfigure}{0.9\textwidth}
\includegraphics[width=\linewidth]{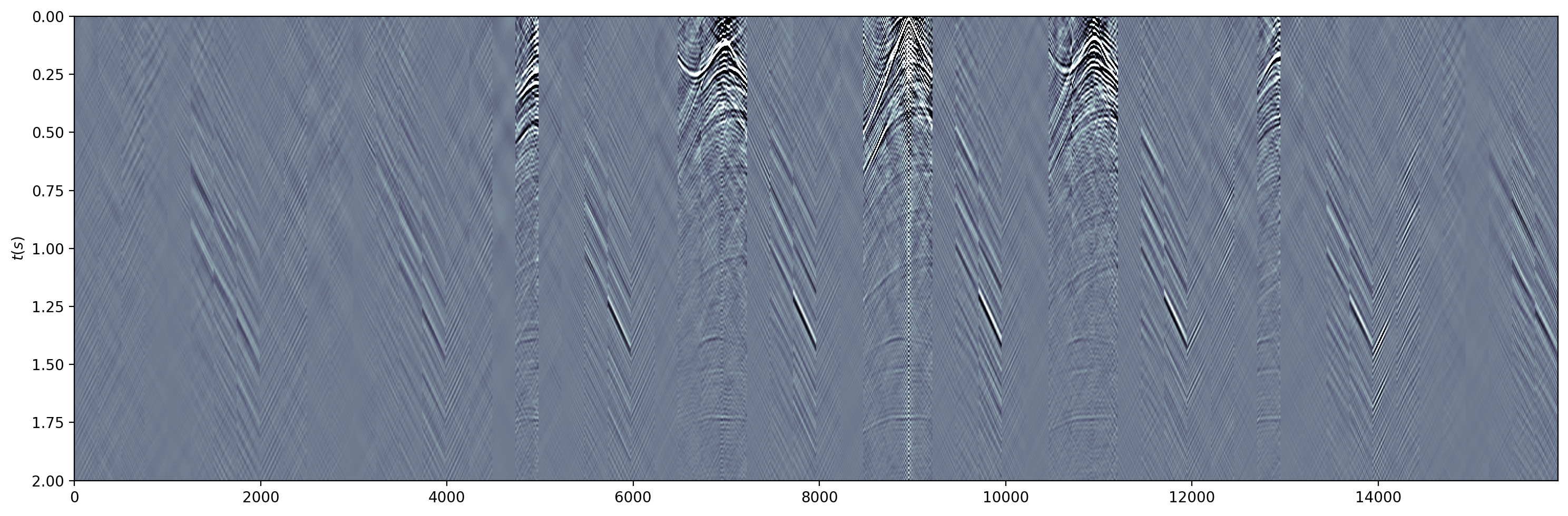}
\caption{$\mathbf{R_{inv}}$, \HT, $r=1$} \label{fig:ds3_r1}
\end{subfigure}\hspace*{\fill}

\caption{MDD results obtained with \HT method for various ranks.} 
\label{fig:ds3}
\end{figure}

Figure~\ref{fig:ds3_r} presents a reference MDD solutions, and Figure~\ref{fig:ds3} illustrates the visualization of inversion results achieved using frequency MDD, applied through method \HT with various rank parameters, compared to the reference results: the modeled solution and the frequency domain full-matrices solution.
This dataset, a substantial 3D one, does not exhibit a globally low rank across its three matrices at all frequencies. Consequently, the only feasible method for this dataset is the \HT method. Our analysis of the previous two datasets shows that the \HT method incurs significant overhead compared to the LR and USV methods. However, it remains preferable over using full matrices. 

\begin{figure}[H] 
\begin{subfigure}{0.34\textwidth}
\includegraphics[width=\linewidth]{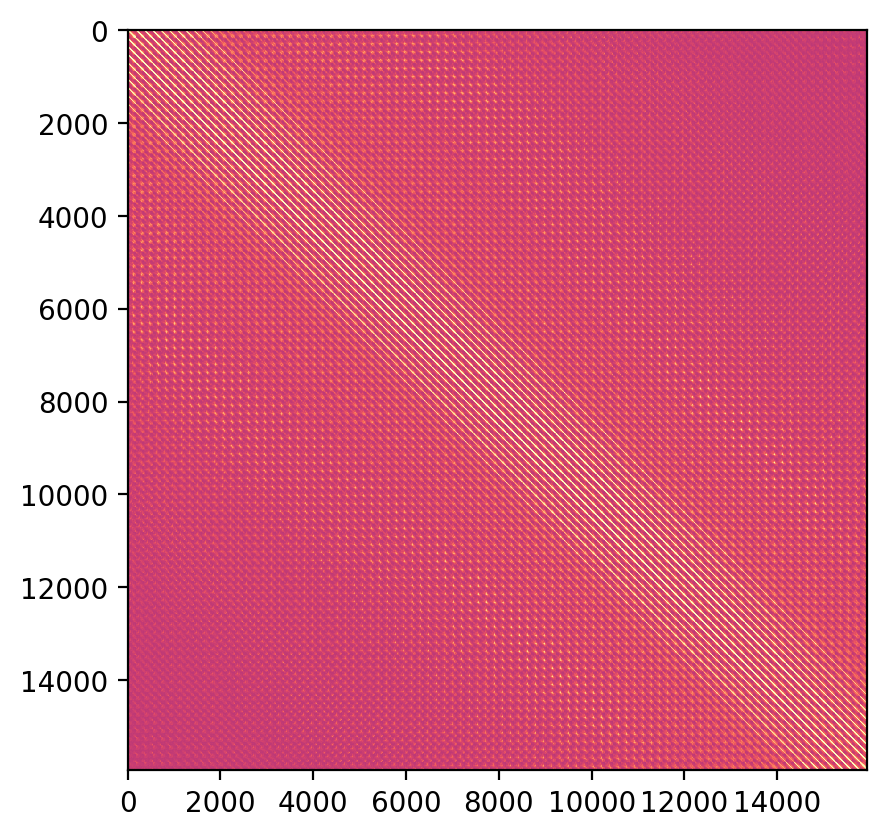}
\caption{Full freq. dom.} 
\end{subfigure}\hspace*{\fill}
\begin{subfigure}{0.34\textwidth}
\includegraphics[width=\linewidth]{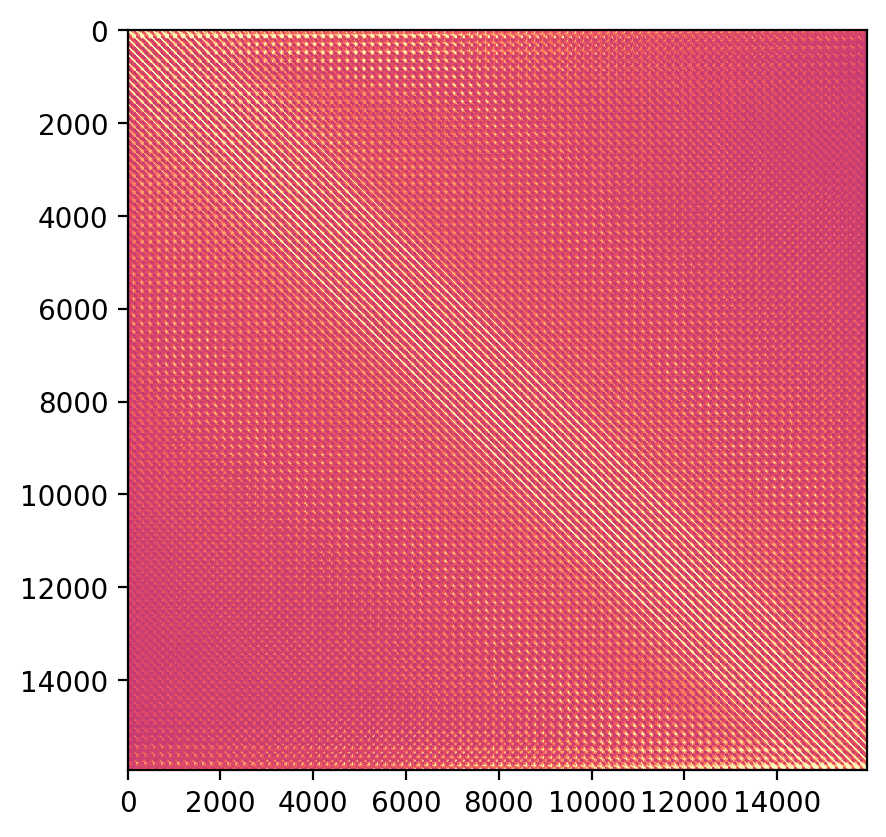}
\caption{\HT, $r=200$}
\end{subfigure}\hspace*{\fill}
\begin{subfigure}{0.34\textwidth}
\includegraphics[width=\linewidth]{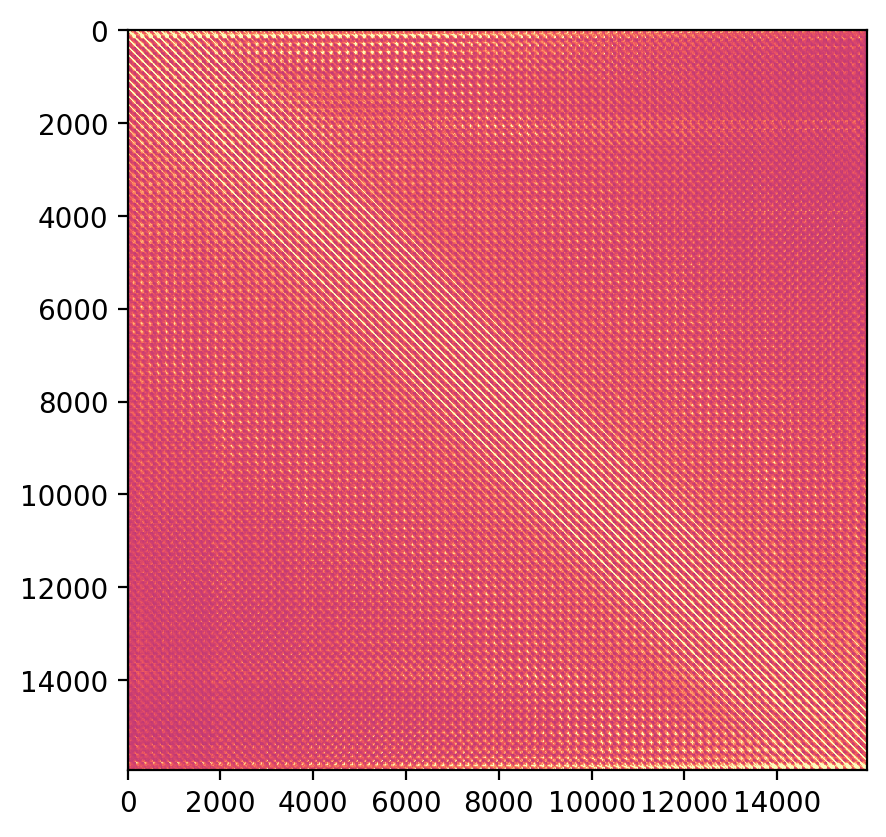}
\caption{\HT, $r=100$} 
\end{subfigure}
\medskip
\centering
\begin{subfigure}{0.34\textwidth}
\includegraphics[width=\linewidth]{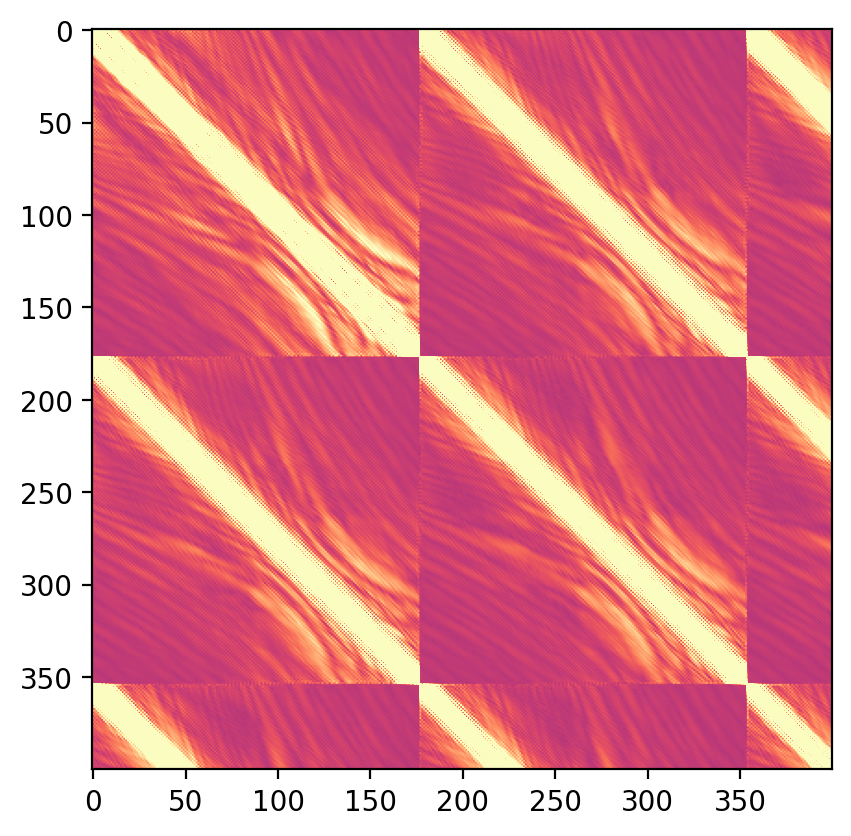}
\caption{Full freq. dom., fragment } 
\end{subfigure}\hspace*{\fill}
\begin{subfigure}{0.34\textwidth}
\includegraphics[width=\linewidth]{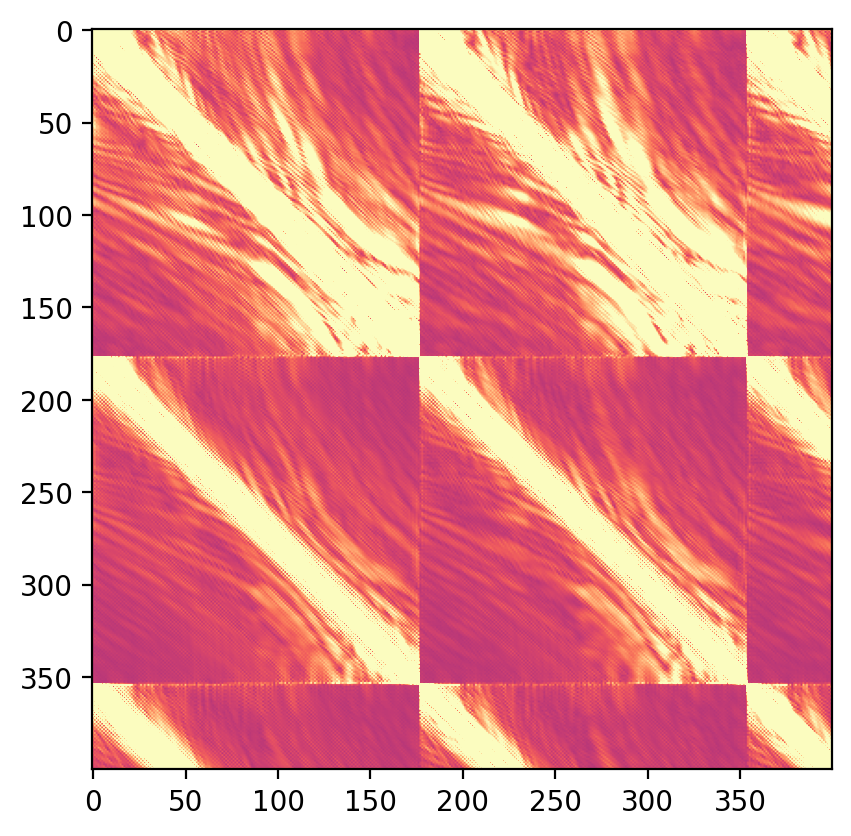}
\caption{\HT, $r=200$, fragment}
\end{subfigure}\hspace*{\fill}
\begin{subfigure}{0.34\textwidth}
\includegraphics[width=\linewidth]{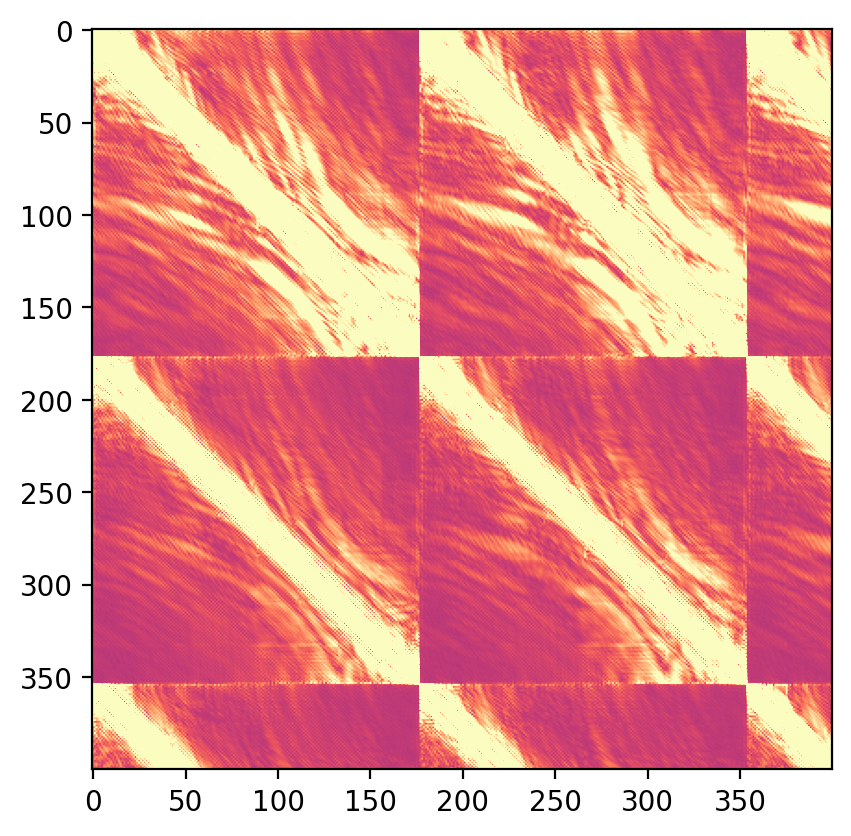}
\caption{\HT, $r=100$, fragment} 
\end{subfigure}
\caption{MDD solution $X_{p}$ in the frequency domain for frequency slice number $p=90$.} 
\label{fig:x_f}
\end{figure}

In Figure~\ref{fig:x_f}, we demonstrate a comparison between the solutions obtained through the full-matrix method and the \HT method within the frequency domain. The frequency slice number is $p=90$. 

This comparison reveals that although the \HT method may represent a compromise in terms of outright solution fidelity, it notably succeeds in maintaining all essential elements of the solution.

Table~\ref{tab:ds3_mem} compares the time and memory requirements for the MDD methods described above. We consider the uncompressed time and frequency domain MDD and \HT method for various ranks. Given the computational limitations of fully resolving the time domain MDD on even the most advanced supercomputers, we are confined to presenting only theoretically derived memory requirements. Conversely, in the frequency domain, leveraging a state-of-the-art supercomputer equipped with parallel GPU cores has enabled us to calculate and visually present the results. However, it is important to note that comparisons of computation time may not be equitable, as the \HT method calculations were performed on a CPU-core-based supercomputer.

For timing, we present integral solution timings, including the approximation to the \HT format. For memory, we sum up the sizes of all matrices involved in the MDD process (operator, right-hand side, and unknown matrices) in the frequency domain. We compute the memory requirements using the formulas~(\ref{eq:mem_rhs_t} --
\ref{eq:mem_x_usv}). The sizes of the \HT unknowns side data are computed empirically. In this 3d Ocean-bottom cable dataset, the number of sources $N_s = 26040$, number of receivers and virtual receivers $N_0 = 15930$, number of time samples $N_t = 1126$, and number of frequency samples $N_f = 200$.

\begin{table}[H]
\centering
\begin{tabular}{c||c|c|c|c|cc}
  & full, time dom.  & full, freq. dom. &  \HT, $r=200$ & \HT, $r=100$ &\HT, $r=1$ \\ \hline \hline
Time (s)& - & - & 4189 & 1530 & 350\\
Mem (TB)            &       3.85       &   2.2            &    0.46      &  0.11        & 0.017       \\
\end{tabular}
\caption{Time and memory comparison of the presented methods.}
\label{tab:ds3_mem}
\end{table}

The trade-off becomes apparent: reducing the rank compromises solution quality (as illustrated in Figures~\ref{fig:ds3} and~\ref{fig:ds3_r}), but it significantly saves time and memory resources.

\section{Conclusion}
To advance the solution of linear systems that involve massive right-hand sides and unknowns, this paper addresses the challenge of solving systems with a data-heavy right-hand side and solution, using the example of Multidimensional Deconvolution (MDD) inversion. We introduce and compare several approaches to compress all matrices within the linear system, significantly improving the solution efficiency, including algorithms based on global low-rank (USV,  LR) and block low-rank ($\mathcal{H}^2$). 
\section{Acknowledgements}
Authors are thankful to the Supercomputing Laboratory and the Extreme Computing Research Center at King Abdullah University of Science and Technology for their support and computing resources.

\section*{References}
\bibliographystyle{abbrv}
\bibliography{lib}

\end{document}